\documentclass[12pt]{article}

\usepackage{amsmath,amsthm,amssymb}

\theoremstyle{plain}
\newtheorem{thm}{Theorem}
\newtheorem{prop}{Proposition}[section]

\newtheorem{lem}{Lemma}[section]

\def\R{\mathbb{R}}

\def\v2{\vskip2mm}
\def\n{\noindent}

\def\({(\!(}
\def\){)\!)}
\def\R{{\bf R}}

\def\a{\alpha}
\def\b{\beta}
\def\e{\varepsilon}
\def\de{\delta}
\def\ga{\gamma}
\def\k{\kappa}
\def\la{\lambda}

\def\th{\theta}

\def\Ga{\Gamma}
\def\La{\Lambda}
\def\Om{\Omega}

\def\pf{{\it Proof.}}
\def\v2{\vskip2mm}
\def\n{\noindent}
\def\z{{\bf z}}
\def\x{{\bf x}}
\def\y{{\bf y}}

\def\be{{\bf e}}
\def\0{{\bf 0}}
\def\pr{{\rm pr}}

\def\tst12{{\textstyle \frac12}}

\def\quadd{\qquad\qquad}
\def\n{\noindent}
\def\beq{\begin{eqnarray*}}
\def\eeq{\end{eqnarray*}}

\def\beqn{\begin{equation}}
\def\eeqn{\end{equation}}

\begin{document}

\begin{center}
{\bf The   transition density  of Brownian  motion  killed on  a bounded set} \\
\vskip4mm
{K\^ohei UCHIYAMA} \\
\vskip2mm
{Department of Mathematics, Tokyo Institute of Technology} \\
{Oh-okayama, Meguro Tokyo 152-8551\\
e-mail: \,uchiyama@math.titech.ac.jp}
\end{center}

\vskip8mm
\n

\vskip2mm
\n
{\it AMS Subject classification (2010)}: Primary 60J65;  Secondary   35K20
\vskip2mm
\n
{\it key words}: heat kernel; 
  exterior domain; transition probability

\vskip6mm

\begin{abstract}
We study the   transition density of a standard two-dimensional Brownian motion killed when hitting a bounded Borel set $A$.  We derive the  asymptotic form  of  the density, say  $p^A_t({\bf x},{\bf y})$,  for large times $t$  and for ${\bf x}$ and ${\bf y}$ in the exterior of $A$ valid uniformly  under the constraint  $|{\bf x}|\vee |{\bf y}| =O(t)$.  Within  the  parabolic   regime $|{\bf x}|\vee |{\bf y}| = O(\sqrt t)$ in particular 
$p^A_t({\bf x},{\bf y})$ is shown to behave like $4e_A({\bf x})e_A({\bf y}) (\lg t)^{-2} p_t({\bf y}-{\bf x})$ for large  $t$, where $p_t({\bf y}-{\bf x})$ is the transition kernel of the Brownian motion (without killing)  and $e_A$ is the Green function for the \lq exterior of $A$' with a pole at infinity normalized so that $e_A({\bf x}) \sim \lg |{\bf x}|$. We also provide  fairly accurate   upper and lower bounds of $p^A_t({\bf x},{\bf y})$ for the case
 $|{\bf x}|\vee |{\bf y}|>t$ as well as  corresponding  results for the  higher dimensions.   
 \end{abstract}
\vskip6mm

\section { Introduction and main results}

Let $A\subset \R^d$, $d\geq 2$ be a bounded and non-polar Borel set and  $A^r$  the set of all regular points of $A$. Denote  by $\Om_A$  the unbounded (fine) component of $\R^d\setminus A^r$. Let $p^A_t(\x,\y)$, $t>0,  \x,\y \in \Om_A$ be the 
transition density of a $d$-dimensional standard Brownian motion in $\Om_A$ killed when it hits $A$. In this paper  we obtain an exact asymptotic form of $p^A_t(\x,\y)$ as $t\to \infty$ that holds true  uniformly   for $\x$, $\y$ subject to the constraint
$|\x|, |\y| =O(t)$. It is shown by   Collet, Martinez and Martin \cite {CMM}  (in which $A$ is assumed to be compact) that as $t\to\infty$
\beqn\label{1}
 p^A_t(\x,\y) \sim \frac{2u(\x)u(\y)}{\pi(\lg t)^2t}  \quad \mbox{if}\quad  d=2 
 \eeqn
and
\beqn\label{2}
p^A_t(\x,\y) \sim \frac{u(\x)u(\y)}{(2\pi t)^{d/2}}\quad \mbox{if}\quad  d\geq 3
\eeqn
 uniformly for  $\x, \y$ in any compact set of $\Om_A$.
Here $u(\x)$ is the unique  harmonic function in  $\Om_A$  that  tends to zero as $\x$ approaches  $A^r$ and is asymptotically equivalent  to $\lg |\x|$  or 1 as $|\x|\to \infty$ according as   $d=2$ or $d\geq 3$; and the symbol  $\sim$ designates the asymptotic equivalence: the ratio  approaches unity in any  specified limit.
On the other hand it is readily verified   that if  $\x$ and $\y$ are  both sufficiently far from  the origin  (depending on $t$ when $d=2$) and not located   in the opposite direction relative to $A$ (in a loose sense), then
\beqn\label{3}
 p^A_t(\x,\y) \sim  p^{(d)}_t(|\y-\x|),
 \eeqn 
 where  $|\x|$ denotes the  Euclidian length of  $\x\in \R^d$  and 
$$p^{(d)}_t(x) = \frac1{(2\pi t)^{d/2}}e^{-x^2/2t} \qquad(x\geq 0).$$
In the case $d\geq 3$  formula (\ref{2})  would be naturally extended to (\ref{3})  for unbounded  $\x$ and $\y$ since $u(\x)\to 1$ as $|\x| \to\infty$ (although there still   remains a problem, the    validity of (\ref{3}) being not obvious when $\x$ and $\y$  in the opposite direction), whereas in the case $d=2$ there is a serious gap between (\ref{1}) and (\ref{3}) and  it is not clear how they are linked together.
Our result, stated below,  extends   formula (\ref{1})  by allowing  $\x$ or $\y$ to become indefinitely large along with $t$ and   makes precise a regime of $t, \x, \y$  for which  (\ref{3}) is valid, and thereby fill the gap between (1)  and (3). The asymptotic form in the parabolic regime $|\x|=O(\sqrt t\,)$ and/or $|\y|= O(\sqrt t\,)$ would  particularly be both interesting and important for an obvious reason: our killed Brownian motion  starting at any fixed point  outside $A$ must be located, at a large time $t$,  inside the annulus  $M^{-1}\sqrt t < |\x| <M\sqrt t$ 
with a high probability  if  the constant $M$ is large.

  The present investigation is motivated by  the study of heat equation in $\Om_A$ in which  a fundamental role is played by   the caloric measure  (the harmonic measure of the heat operator $\frac12 \Delta - (\partial/\partial t)$ in the space-time cylinder  $\Om_A\times [0,\infty)$),  which  consists of two parts: one is  the transition probability $p^A_t(\x,\y)|d\y|$,  the part on the  initial boundary $\Om_A\times \{0\}$  and  the other,  the lateral part, is   
the hitting distribution in space-time of the boundary $\partial A$. Here $|\cdot|$ designates the Lebesgue measure on $\R^d$ and $\partial A$  the  Euclidean boundary of $A$.
 The asymptotic form of the lateral part is computed in \cite{Ucalm}  and  our present result together with it identifies the explicit asymptotic form of the caloric measure  for large time valid uniformly at least in the  regime  $|\x|\vee |\y| = o( t)$, where $x\vee y = \max\{x,y\}$. 

\v2
  Let  $U(a)\subset \R^d$ denote  the $d$-dimensional open ball  centered at the origin with radius $a>0$.  $ P_\x$ denotes the law of  a Brownian motion $(B_t)$  started at $\x$.  For a bounded Borel set $A$,   $\sigma_A $ 
denotes the first hitting time of $A$ by $B_t$, namely  $\sigma_A= \inf\{t>0: B_t\in A\}$, and  $A^r$ the set of all regular points of $A$, namely  the set of  $\y$ such that
$P_\y[\sigma_A=0]=1$. Put $R_A = \sup \{|\z|: \z\in A^r\}$.  The set $\Om_A$ is then defined by   
$$\Om_A = \{\x\in \R^d : P_\x[\sigma_{\partial U(R)}< \sigma_A] >0 \} $$
with   any/some $R>R_A$, and the transition density $p^A_t(\x,\y)$, $\x, \y \in \Om_A$ by
$$p^A_t(\x,\y) = P_\x[B_t\in d\y, \sigma_A>t]/|d\y|.$$  
   $p^A_t(\x,\y)$ is symmetric  and positive for all $\x, \y\in \Om_A$  and $t>0$ and jointly continuous  in $(t,\x,\y)\in (0,\infty)\times \Om_A^\circ\times \Om_A^\circ$ ($\Om_A^\circ$ denotes the interior of $\Om_A$).
 Readers may refer to \cite{B} (Section 2.4, especially Theorem 4.4 and Proposition 4.9) for  existence  and  properties of $p^A_t(\x,\y)$.  
The set  $A^r$  is Borel 
 and   $A\setminus A^r$ is polar, so that  $P_\x[\sigma_A =\sigma_{A^r}]=1$ and  $P_\x[B_{\sigma(A)}\in A^r|\sigma_A<\infty]=1$ for all $\x$ (see e.g. \cite{B}, \cite{BG}, \cite{C}).
Denote   by ${\rm nbd}_\e (A^r)$ the open    $\e$-neighborhood of  $A^r$ in $\R^d$.  
We also suppose
$$ \R^d \setminus \Om_A = A^r$$
($A^r$ has no \lq fine cavity'), which requires no  essential change of the content of the paper
since we are concerned only the exterior problem. Because of it  we may write $\x\notin {\rm nbd}_\e (A^r)$
 for   $\x \in \Om_A\setminus {\rm nbd}_\e (A^r)$. We  write $x$ for the Euclidean length of $\x\in\R^d$: $x=|\x|$.

Let $d=2$ and 
$e_A(\x)$ be the Green function for $\Om_A$ with a pole at infinity normalized so that  $e_{A}(\x) \sim \lg x$ when  $ x \to \infty$ (see (\ref{8}) for a definition). The following theorem is stated only  in the case $x\leq y$  which  is not a restriction because  $p^A_t(\x,\y)$ is symmetric in $\x$ and $\y$ as mentioned above.

\begin{thm} \label{thm1}\, Let $d=2$.  For any $M\geq 1$ and  $\e>0$,  uniformly for $\x,\y \notin {\rm nbd}_\e (A^r)$ with  $x \leq y<Mt$, as $t\to \infty$
\beqn\label{4}
\frac {p^A_t(\x,\y)}{p^{(2)}_t(\y-\x)}  \sim \left\{ \begin{array} {ll}
\displaystyle{ \frac{4e_A(\x)e_A(\y)}{(\lg t)^2}} \quad & \mbox{if} \quad y\leq M\sqrt t, \\[4mm]
\displaystyle{ \frac{e_A(\x)}{\lg (t/y)}} \quad & \mbox{if} \quad y> \sqrt t,  \, xy \leq Mt \,\,  \mbox{and}\,\, y/t\to 0, \\[4mm]
1 \quad & \mbox{if} \quad  \,  xy \geq t / M\,\, \mbox{and}\,\, x \to \infty. 
\end{array} \right.
\eeqn
\end{thm}

\v2\v2
The function $e_A$ agrees with $u$ appearing  in (\ref{1})  (verification is standard and given  in \cite{CMM}), so that  Theorem \ref{thm1}  improves  (\ref{1}) by extending the range of validity to the regime  $x\vee y<M\sqrt n$.  The condition $\x,\y \notin {\rm nbd}_\e (A^r)$ can be relaxed to $\x,\y \in \Om_A$ under some regularity condition for $\partial A$ (smoothness is sufficient: cf. Theorem \ref{thm2.1} in Section 2.1 and a comment right after it), but not in general. 

The case when $y > \de t$ and $x< 1/\de$  for a constant $\de>0$ is excluded from the condition  $x \leq  y < Mt$ in Theorem \ref{thm1}.  In this case   relation (\ref{4}) breaks down (with the condition $y/t\to0$ violated in the second case or  the condition $x\to\infty$  in the third): in fact   the second and third relations in 
(\ref{4})  are complemented by the following result: there is a continuous function   
$c^A(\x;{\bf v})$ taking values in $(0,1)$ such that  under the constraints  $y \asymp t$ and $x \asymp 1$, as $t\to \infty$
$${p^A_t(\x,\y)}\; \sim \;c^A(\x; \y/t)\,{p^{(2)}_t(\y-\x)}$$
(Proposition \ref{prop2.1}), which combined with Theorem  \ref{thm1} shows
that  if $\x,\y \notin {\rm nbd}_\e (A^r)$ (for an $\e>0$) and if 
$$ \mbox{ $ x\vee y < Mt$ \; and \; $1<t < M (x+1)(y+1) $}$$
for a constant  $M\geq 1$,  then
\beqn\label{eq1.2} {p^A_t(\x,\y)}\asymp {p^{(2)}_t(\y-\x)}.
\eeqn
Here the expression $f(t)\asymp g(t)$  means that the ratio $f(t)/ g(t)$ is bounded away from both zero and infinity (provided that $f(x) g(x) > 0$). Relation  (\ref{eq1.2}) is a lower bound since the upper bound is trivial and virtually follows if we show it with $A$ replaced by a disc $U(a)$ that includes 
$A$.

Let $\be=(1,0)$, the first coordinate vector of $\R^2$.    The next theorem provides upper and lower bounds  when  $A=U(a)$ and $x\vee y>\!\!>  t$. (Here and in the sequel  the expression  $a<\!\!< b$ ($a>0$), used in  informal expository passages, means that $b/a$ is  large enough.) We have only to consider 
$\y=- y\be.$ 
 Write $\x=(x_1, x_2)$. It being readily shown (see proof of Lemma \ref{lem3.2}) that  if
  $A\subset \overline{U(a)}$, then  (\ref{eq1.2}) is true  under the condition 
$$
a^2< t< \de ay,\; x> (1+\de)a    \mbox{\; (for some $\de>0$)\;
and \; either $x_1\leq 0$ or  $|x_2|\geq a$,}
$$
we concentrate on the case when  $x_1>0$ and $|x_2|<a$. For $x_1>0$, define the function  $\rho=\rho(\x,\y)$ by
$$\rho = \frac{x_1}{x_1+y}.$$ 

\begin{thm} \label{thm2}\,  Let $d=2$. Suppose  $t >a^2$,  $x \leq y$,  $\y=- y\be$, $x_1>0$ and $0\leq x_2<a$.  
\v2
{\rm (i)} \,  For each  $\de>0$, (\ref{eq1.2}) holds if  $A\subset \overline{U(a)}$, $\x \notin U(a(1+\de))$ and 
$$ \rho t\geq  \de a^2\wedge (a- x_2)^2;$$

{\rm (ii)} \, For any $\de>0$ there exist  a constant  $\k_{\de}$ depending only on $\de$ and a universal constant $c>0$  such that  if $\x \notin U(a+\de\sqrt{\rho t}\,)$ and $ a- x_2 \geq \de  \sqrt{\rho t}$, then 
$$\frac{\k_\de\,\sqrt {\rho t}}{a- x_2}\exp\bigg\{-k_*(\x)\bigg(\frac{a- x_2}{\sqrt{\rho t}} + \de\bigg)^2\bigg\} 
\leq \frac{p^{U(a)}_t(\x,\y)}{p^{(2)}_t(\x-\y)} \leq      \frac{c\,\sqrt {\rho t}}{a- x_2}\exp\bigg\{-k(\x)\frac{(a- x_2)^2}{\rho t}\bigg\}
$$
with some  functions $k_*(\x)$ and  $k(\x)$  that satisfy     
$$\sqrt2- 1< k(\x) < \tst12 <k_*(\x)<1 \quad{and}\quad  \lim_{x_1\to\infty} (\tst12-  k(\x))x_1^2 = \lim_{x_1\to\infty}  (k_*(\x) - \tst12)x_1^2 <\infty.$$ 

\v2
{\rm (iii)} \, If  $x_1>2a $ and $ a^2<t < ay$, then 
$$ \frac{p^{U(a)}_t(\x,\y)}{p^{(2)}_t(\x-\y)} \asymp  \bigg[\frac{\sqrt {\rho t}}{a- x_2}\wedge 1\bigg]
\exp\bigg\{-\Big(\frac12 + O(a/x_1)\Big)\frac{(a- x_2)^2}{\rho t}\bigg\} $$
where the  constants involved in $\asymp$ as well as in  the $O$ term  are universal. 
\end{thm}

\v2
In the case when  $ (a-x_2)^2/\rho t \to\infty$ and $a/x_1 \to 0$,  $(2\pi)^{-1/2}$ times 
the right side of the  formula in (iii) above
may be  heuristically regarded as an asymptotic form of   the probability that  a Brownian bridge   joining $\x$ and $\y$ with  length $t$ crosses   the vertical axis at a level higher than  $a$,
for on the one hand if $B'_t$ designates the second component of $B_t$,
\beqn\label{505}
\sqrt{2\pi} P_\x[B'_{\rho t}>a ] \sim \frac{\sqrt{\rho t}}{a-x_2} e^{- (a-x_2)^2/{2\rho t}},
\eeqn
  and  on the other hand  the distribution of  the crossing  time  almost   concentrates at $\rho t$.

Higher dimensional analogues of Theorems \ref{thm1} and \ref{thm2} will be given in Sections 2.4 and 3.3, respectively. The proofs are almost the same as to the two dimensional results. However,  the formula corresponding to  (\ref{505}) is not similar to it and the result corresponding to Theorem \ref{thm2} is accordingly modified.
\v2\n
{\sc Remark 1.} \, In the regime $x\vee y =O(t)$ the asymptotic behavior of $p_t^A(\x,\y)$ is almost the same as that of  $p_t^{U(a)}(\x,\y)$ if  $x\wedge y$ is large enough  so that $e_A(\x)$ and  $e_A(\y)$  
are well approximated by $\lg x$ and $\lg y$, respectively. The situation becomes quite different when 
$(x\vee y)/t \to \infty$: e.g., if  the length of ${\rm pr}_{\be} A$  (the projection to the line perpendicular to $\be$)  is zero and $\partial A$ is smooth, then $p_t^A(\x,y\be)\sim p^{(2)}_t(y\be -\x)$  
 as $y/t\to \infty$ for all $\x \in \Om_A\setminus  {\rm nbd}_\e(A^r)$ (cf. \cite[Theorem 3.5]{Ucalm}). 
 
 \v2

 In \cite{CMM}  Collet, Martinez and Martin   employ  a  Harnack inequality and a classical expression of  $p^{U(1)}_t$ as  given in \cite{CJ} to obtain an upper bound of  $p_t^A$ and then  use a compactness argument.  The limiting  function is identified by the uniqueness of the harmonic function $u$.

 
 Our approach is quite different.   If  restricted to the first case  of (\ref{4}), our proof  is based on the upper bounds  
\beqn\label{000}
\frac{P_\x[\sigma_{U(1)}\in dt, B_t \in d\xi]}{dt |d\xi|} \leq \frac{C}{t  \lg t} \qquad (1< x <  \sqrt t\, , d\xi\subset \partial U(1)),
\eeqn
where $|d\xi|$ designates the length of $d \xi$, and
\beqn\label{5}
 P_\x[\sigma_{U(1)} <t]  = o(1) \quad \mbox{if (and only if)} \quad  \liminf \frac{\lg x}{\lg t} \geq \frac12 
\eeqn
(as $t, x \to\infty$)
 as well as on
the following relation: for each  $\e>0$, uniformly for $\x\in U(r)\setminus $nbd$_\e(A^r)$,  
 \beqn\label{6}
P_\x[\sigma_{\partial U(r)} < \sigma_A] \;  \sim \;  \frac{e_A(\x)}{\lg r}  \quad
 \mbox{ as} \quad  r\to \infty.
    \eeqn  
Relations  (\ref{000}), (\ref{5}) and  (\ref{6})  are taken  from     \cite{Ubh}, \cite{Ucal_b} and \cite{Ucalm}, respectively. 
Roughly speaking we first obtain (\ref{4}),  by using  (\ref{000}) and (\ref{5}),     under the condition $x\sim y \sim \sqrt {t}/ \lg t$ (a special case of the first one of (\ref{4})). This means 
that the  Brownian bridge  that connects   $\x$ and   $\y$  visits $A$ with so small a probability that 
$p^A_t(\x,\y) \sim p_t^{(2)}(0)$ in this special case. 
Next we consider the case  $x<y\sim \sqrt t/\lg t$ and the event that  the process  starting at $\x$ escapes from  $A$ before exiting  the disc  $U(\sqrt t/\lg t)$ with the exiting time $ o(t)$. Formula  (\ref{6})
says that the probability of this escape from $A$ is asymptotic to $e_A(\x)/ \lg \sqrt{t}$, which together with the preceding result shows
 $p^A_t(\x,\y)\sim 2(\lg t)^{-1}e_A(\x)p_t^{(2)}(0)$. 
For the case  $x\vee y <\!\! <\sqrt t/\lg t$, we have only to apply the latter case result to the  Brownian motion started afresh at the exit  time from the disc $U(\sqrt t/\lg t)$. These are to   verify  (\ref{4}) in the case  $x\leq y \leq C \sqrt t /\lg t$.

 The proofs for the other cases of  (\ref{4}) need some  results from \cite{Ucal_b}  and \cite{Ucalm} in addition to those mentioned above.  The same  method as described above   applies to  random walks as is discussed in  \cite{Udnst}.

The
ratio on the left side of (\ref{4}), $p_t^A(\x,\y)/p_t^{(d)}(\y-\x)$,  is  the probability that the Brownian bridge joining  $\x$ and $\y$ in the time interval $[0,t]$ avoids $A$, namely
$$P_\x[ B_s\notin A  \,\, \mbox{for}\,\, 0\leq s \leq t\,|\, B_t= \y].$$
Thus  formula (\ref{4}) may be considered as giving the asymptotic form of this probability. 
 In view of Theorem A.2 of Section 4 we have $P_\x[\sigma_A > t]  \sim 2e_A(\x)/ \lg t$ if  $\limsup\frac{ \lg x}{ \lg t} \leq \frac12$. If $\limsup\frac{ \lg (x\vee y) }{ \lg t} \leq \frac12$,     the right side  of (\ref{4}) is therefore  asymptotic to 
 \beqn\label{7}
 P_\x[\sigma_A >t]P_\y[\sigma_A>t],
 \eeqn
 the probability that two independent Brownian motions starting at $\x$ and $\y$ both avoid $A$ until  the time $t$.  Collet et al.  \cite{CMM} points out  (for their case)   that (\ref{1}) is interpreted as stating the coincidence of the asymptotic forms  of    these two probabilities  for  $\x$ and $ \y$ restricted to 
  a compact set of $\Om_A$. 
 Formula  (\ref{4})  shows  that  this coincidence extends over  the regime  in which 
 $\limsup\frac{ \lg (x\vee y) }{ \lg t} \leq \frac12$ but not necessarily beyond  it.   It may be worth noting that  since the law of 
 the Brownian Bridge is not affected by a constant drift, for any constant vector ${\bf v}$ the transition density of the process $B_s^{{\bf v}} := B_s+{\bf v}s$ killed upon hitting  $A$ is given by
 $$
  \frac {p^A_t(\x,\y)}{p^{(d)}_t(\y-\x)}  p^{(d)}_t(\y-{\bf v}t-\x),$$
  hence Theorem \ref{thm1} also determines its asymptotic form.

 We separate   Theorem \ref{thm1}  into three propositions (Theorem 3, Proposition 2.1 and 2.2) corresponding to  the three cases  in (\ref{4}), namely  cases (1) $x\vee y \leq \sqrt t$, (2)  $xy \leq t$, $x\vee y > \sqrt t$ and (3)  $xy> t$,  and deal with them in Sections 2.1, 2.2 and 2.3, respectively. 
 The result  corresponding to Theorem \ref{thm1} in  the higher dimensions $d\geq 3$ will be briefly  discussed in Section 2.4,  the proof being  much simpler.   Theorem \ref{thm2} is proved in Section 3. The higher dimensional result for it is given in Section 3.4; the result itself  somewhat differs from the two dimensional one although its proof is  similar.
  In Section 4 we collect the results from \cite{Ubh}, \cite{Ubes}, \cite{Ucal_b} and \cite{Ucalm} that are applied in this paper.

 \section{Asymptotic forms in the regime $x\vee y =O(t)$}

 For $\x\in \Om_A$, define
\beqn\label{H_A}
H_A(\x,t; d\xi) = \frac{P_{\x}[ B_{\sigma_A}\in d\xi, \sigma_A\in dt]}{dt}\quad\quad (d\xi \subset \partial A, t>0),
\eeqn
the hitting distribution of $A$ in space-time, which  in addition to the things
introduced in Section 1  will be of fundamental use throughout  the rest of the paper.  We shall apply
some asymptotic estimates of  $H_A$ (or rather its partial integrals) that are taken from \cite{Ubh}, \cite{Ubes}, \cite{Ucal_b} and \cite{Ucalm}.  We collect them in Section 4 for convenience of citation.  
 
 Let $d=2$. The function  
$e_A(\x)$, $\x \in \Om_A$ may be defined  by   
   \beqn\label{8}
   e_A(\x)= \pi \lim_{|\y|\to\infty} g_{\Om_A}(\x,\y),
   \eeqn 
where   
$g_{\Om_A}(\x,\y) = \int_0^\infty p^A_t(\x,\y)dt\,$ $  (\x,\y \in \Om_A)$, the Green function for 
the set  $\Om_A$;
  in particular $e_A$ is harmonic in the interior of  $\Om_A$.  Let ${\rm err}_*(\x,t)$ stand for any
  function  satisfying that  ${\rm err}_*(\x,t)=0$  if $x:=|\x|\geq 2R_A$ and  
  $$|{\rm err}_*(\x,t)| \leq C P_\x\Big[\sigma_{A\cup \partial U(2R_A)} > t/\lg t\,\Big]\quad\mbox{if}\quad x < 2R_A.$$
   It is noted that the right side above is dominated by $Ce^{-\la t/R_A^2\lg t}$ with a certain  universal constant $\la>0$ ($*$ is attached to distinguish from the notation ${\rm err}(\x,t)$  used in \cite{Ucalm}).  
  
   We  often write  $p^{(d)}_t(\x)$  for  $p^{(d)}_t(x)$  for convenience sake; also write $\sigma(A)$ for $\sigma_A$ for typographical reason; $x\vee y$ and $x\wedge y$ denote  the maximum and  minimum of real numbers $x, y$, respectively.  

\subsection{ Case  $x\vee y =O(\sqrt t\,)$, $d=2$}

 In this and the next subsections we assume  $d=2$  and  simply write  $p_t(x)$  for  $p_t^{(2)}(x)$ (and also  $p_t(\x)$  for  $p_t^{(2)}(\x)$).
The following theorem refines the first case of (\ref{4}).     \v2
\begin{thm} \label{thm2.1}\,  For any $M>1$,  uniformly for $\x,\y \in \Om_A$   
subject to the constraint $x\vee y < M \sqrt{t}$,   
as $t\to \infty$
\beqn\label{9}
p^A_t(\x,\y) 
=\frac{4e_A(\x)e_A(\y)}{(\lg t)^2}p_t(\y-\x) \bigg(1+ O\bigg(\frac{\lg\lg t}{\lg t}\bigg) \bigg) + {\rm err}_*(\x,t) +{\rm err}_*(\y,t). 
\eeqn
\end{thm}

The constants involved in the $O$ term appearing in (\ref{9}) of course depends on $M$ and the same remark applies to the results that follow. 
 If  the boundary $\partial A$ is smooth, then ${\rm err}_*(\x,t)\leq C e_A(\x) e^{-\la\sqrt t/ R_A^2\lg t}$ 
so that ${\rm err}_*(\x,t) +{\rm err}_*(\y,t)$ may be deleted from the right side of  (\ref{9}) unless  $\x$ and $\y$ approach $\partial A$ simultaneously. 

The proof of Theorem \ref{thm2.1}   will be given at the end of this subsection after showing  Lemmas \ref{lem2.2} and \ref{lem2.3} given below.  Lemmas \ref{lem2.2} and \ref{lem2.3} improve the first relation of (\ref{4}) by giving  error estimates  and extending  the range of validity of the formula beyond the restriction $x\vee y = O(\sqrt t)$.  In order to include the improvement   we shall apply certain refined versions of (\ref{000}), (\ref{5}) and (\ref{6});  otherwise    we do not need such refined ones.  We suppose $d=2$ throughout this and the next  subsections.

\begin {lem}\label{lem2.2} For any  $M\geq 1$,  as $x\wedge y\wedge t\to\infty$ under the constraints
$\lg x\sim\lg y\sim \tst12 \lg t$ and  $xy\leq Mt$,
$$p^A_t(\x,\y) = p_t(\y-\x)\bigg[1+ O\bigg(\frac1{\lg t}\Big[1\vee \lg\frac{t}{x^2\wedge y^2}\Big]\bigg)\bigg].$$
\end{lem}
\v2\n
\pf\, Since $p_t^{A}(\x,\y) \leq p_t^{U(R_A)}(\x,\y) \leq p_t(\y-\x)$, it suffices to obtain a proper upper bound of 
the difference $p_t(\y-\x) - p_t^{U(R_A)}(\x,\y) $. To this end we use  the identity 
\beqn\label{10} 
p_t(\y-\x) - p_t^{U(R_A)}(\x,\y) =  \int_0^t ds\int_{\partial U(R_A)} H_{U(R_A)}(\x,s;d\xi)p_{t-s}(\y-\xi).
\eeqn
Let $I_{[a,b]}$ denote the double integral in (\ref{10})  restricted to  $[a,b]\times \partial A, 0\leq a<b\leq t$.
Suppose   $R_A=1$ for simplicity.  Let $xy<Mt$.  We may suppose $x\leq y$, so that  $x<\sqrt{Mt}$.   On noting 
$p_{t-s}(\y-\xi) \leq  Cp_t(y)$ for $\xi \in \partial U(1), s<t/2$,
\[
I_{[0,t/2]}  \leq C p_t(y)P_\x[\sigma_{U(1)}\leq t/2].
\]
  We deduce from Theorem A.2 of Section 4  (with $U(1)$ in place of $A$)  that if $\lg x\sim\frac12 \lg t$
 $$P_\x[\sigma_{U(1)}\leq t/2] \leq   \frac{C[1\vee  \lg (t/x^2)]}{\lg t};$$ 
  we also have 
$p_{t}(y) \leq p_t(\y-\x) e^{2M}$. Putting these together we obtain 
  \beqn\label{101}
I_{[0,t/2]}  \leq C e^{2M} \frac{1\vee \lg (t/x^2)}{\lg t} p_t(\y-\x).
\eeqn
For the other part  $I_{[t/2,t]}$ we  apply  Theorems A.1 and A.4 to see that if $\lg x \sim \lg y \sim \frac12 \lg t$, then
$$I_{[t/2,t]} \leq  \frac{Cp_t(x)}{\lg t}\int_{\partial U(1)}|d\xi|\int_{t/2}^t p_{t-s}(\y-\xi)ds. $$
Writing $r$ for $|\y-\xi|$ and changing the variable by $u=r^2/2(t-s)$ so that $ds/(t-s) =du/u$, we  compute the inner integral to obtain an upper bound of it as follows:
$$\int_{t/2}^tp_{t-s}(r) ds=\frac1{2\pi} \int_{r^2/t}^\infty  \frac{e^{-u}}{u}du 
\leq C\Big[e^{-r^2/t} \vee \lg \frac{t}{r^2} \Big] \leq C' e^{-r^2/t} \Big[1\vee \lg\frac{t}{r^2} \Big],$$
and then   infer that 
\beqn\label{**}
I_{[t/2, t]} \leq  Ce^{M} p_t(\y-\x)\frac{1\vee \lg (t/y^2)}{\lg t}.
\eeqn
Combined with (\ref{101}) this  concludes  the proof  of the lemma. \qed


\v2\v2
\begin {lem}\label{lem2.3} \, For   $M\geq 1$,  if $\x\in \Om_A$, $x < M \sqrt{t}/\lg t$, $y <M\sqrt{t}$ and   $\lg y\sim \frac12 \lg t$, then 
$$p^A_t(\x,\y) = \frac{2 e_A(\x) + {\rm err}_*(\x,t)}{\lg t} p_t(y) (1+o(1));$$
moreover if $ y \geq \sqrt{t}/(\lg t)^\de$  for some  $\de>0$   in addition, the error term $o(1)$ above can be replaced by  
$$O((\lg\lg t)/\lg t).$$
\end{lem}
\v2\n
\pf\, Let  $x< r:=M\sqrt{t}/\lg t$, $y< M\sqrt t$ and   $T:= 2 t/\lg  t$
and   decompose
\begin{eqnarray} \label{13}
p^A_t(\x,\y) &=& \int\!\!\int_{[0,T]\times \partial U(r)} P_\x[ \sigma_{\partial U(r)}\in ds, \sigma_A>s, B_{\sigma_{\partial U(r)}}\in d\xi \,]\, p^A_{t-s}(\xi,\y) \\
&& \, + \, \e(t,\x,\y),\nonumber 
\end{eqnarray} 
where
$$\e(t,\x,\y)=\int_{U(r)\cap \Om_A} P_\x[ \sigma_{\partial U(r)}\wedge \sigma_A >T, B_T\in d\z] \, p^A_{t-T}(\z,\y).$$

 
  It holds that
\beqn\label{131}
 P_\x[\sigma_{\partial U(r)}\wedge \sigma_A >T] \leq   C e^{-\la T/r^2}\Big(e_A(\x) + P[\sigma_{\partial U(2R_A)}\wedge \sigma_A >{\textstyle \frac12} T]\Big)
 \eeqn
with  some universal constant $\la>0$. 
Indeed, this is trivial if $x> R:=2R_A$ for which $e_A(\x)> \lg 2$ (cf.(\ref{asymp_e})), whereas,
by applying strong Markov property at the exit  time  $\sigma_{\partial U(R)}\wedge \sigma_A=\sigma_{A\cup\partial U(R)}$ and breaking the event $\sigma_{\partial U(r)}\wedge \sigma_A >T$ according as 
 $\sigma_{ \partial U(R)}$ is less than $T/2$ or not,  we  infer 
 that for $\x \in U(R)\setminus A^r$, 
\[
P_\x[\sigma_{\partial U(r)}\wedge \sigma_A >T] 
 \leq
 \Big( P_\x[\sigma_{\partial U(R)} < \sigma_A] + P_\x[\sigma_{A\cup\partial U(R)} >{\textstyle \frac12} T] \Big)\sup_{\y\in U(R)}P_\y[\sigma_{\partial U(r)} >{\textstyle \frac12} T].
 \]
 The supremum 
above is $O(e^{-\la T/r^2})$;  use Theorem A.3 to see $P_\x[\sigma_{\partial U(R)} < \sigma_A]\leq e_A(\x)/\lg 2$. Thus  we have (\ref{131}).

By  $y\leq M\sqrt{t}$ we have
\beqn\label{15}
ry/t \leq M^2/\lg t.
\eeqn 
Then,  noting  $p^A_{t-T}(\z,\y)\leq p_{t-T}(\y-\z)\leq c_M p_{t}(y)$ for  $z\leq r$, we observe 
\begin{eqnarray} \label{14} 
\e(t,\x,\y) 
\leq  c_M'e^{-\la T/r^2}[e_A(\x) + {\rm err}_*(\x,t)]p_t(y), 
\end{eqnarray}
hence $\e(t,\x,\y)$ is absorbed into the error terms  since  $T/r^2 \geq M^{-2}\lg t$. 

For $s<T, \, \xi \in \partial U(r)$, by Lemma \ref{lem2.2} we have
\beqn\label{***} 
p^A_{t-s}(\xi,\y) = p_{t-s}(\y-\xi)(1+ o(1));
\eeqn
 observing    $T|\y-\xi|^2/t^2 \leq 2M^2/\lg t$,
from (\ref{15})  we deduce
\beqn\label{p} p_{t-s}(\y-\xi) = p_t(y)\Big(1+O\Big(\frac1{\lg t}\Big)\Big). 
\eeqn
On using Theorem A.3 of Section 4 as well as  (\ref{131})
\begin{eqnarray}\label{24} 
P_\x[\sigma_{\partial U(r)} <T \wedge \sigma_A] &=& P_\x[\sigma_{\partial U(r)} < \sigma_A] - P_\x[T< \sigma_{\partial U(r)}  < \sigma_A] \nonumber\\
&=& \frac{e_A(\x) + {\rm err}_*(\x,t)}{\lg r}\bigg( 1 + O\bigg(\frac1{\lg r}\bigg)\bigg).
\end{eqnarray}
Now the double integral in (\ref{13}) may be written as 
\beqn\label{28}
 P_\x[\sigma_{\partial U(r)} <T \wedge \sigma_A]\, p_t(y) (1+ o(1))
=
 \frac{e_A(\x)  + {\rm err}_*(\x,t)}{\lg r}  p_t(y) (1+ o(1)).
\eeqn
  Since  $\lg r = \frac12 \lg t + O(\lg\lg t)$, this verifies the first half of the lemma.  
 
 If $ y \geq  \sqrt{t}/(\lg t)^\de$, then   in  (\ref{***})    $o(1)$  can be replaced by $O((\lg\lg t)/\lg t)$ according to Lemma \ref{lem2.2} (note by  (\ref{15}) we also have  $|\y\cdot\xi|/t \leq ry/t= O(1/\lg t)$). 
This permits replacing $o(1)$ by $O((\lg\lg t)/\lg t)$ in (\ref{28}), and thereby yields the required error bound.
  \qed
\v2
{\sc Proof of Theorem \ref{thm2.1}.} The case $y\geq x\geq \sqrt t/\lg t$ is covered by Lemma \ref{lem2.2} (see (\ref{asymp_e}) in Section 4 for the error arising from the  replacement of $e_A(\x)$ by $\lg x$).  First suppose $x\leq y\leq M\sqrt{t} /\lg t$. Then, by  Lemma \ref{lem2.3}  we see  in  the double integral in (\ref{13})
 $$p^A_{t-s}(\xi,\y) = p_{t-s}^A(\y,\xi) = \frac{2e_A(\y) + {\rm err}_*(\y,t)}{\lg t} p_{t-s}(\xi) \bigg(1+ O\bigg(\frac{\lg\lg t}{\lg t}\bigg) \bigg). $$
Since  $p_{t-s}(\xi) = p_t(0)(1+O(1/\lg t))$ for $s\leq T$, 
we may write (\ref{13}) as 
$$p^A_t(\x,\y) = P_\x[\sigma_{\partial U(r)} < T \wedge \sigma_A]\frac{2e_A(\y) + {\rm err}_*(\y,t)}{\lg t} p_t(0) \bigg(1+ O\bigg(\frac{\lg\lg t}{\lg t}\bigg) \bigg) + \e(t,\x,\y).  $$
 The evaluation of the term  $\e(t,\x,\y) $ given in (\ref{14})   and that of
  $P_\x[\sigma_{\partial U(r)}  < T\wedge \sigma_A]$   in (\ref{24}) are valid. 
  Thus   we conclude the required relation of the Theorem \ref{thm2.1}. Finally in the case $x\leq \sqrt t /\lg t \leq  y$, we have $xy \leq Mt/\lg t$, so that   $p_t(y) = p_t(\y-\x) + O(1/\lg t)$ and applying the second half of Lemma \ref{lem2.3} finishes the proof. \qed
 
   \subsection{Case  $xy =O(t)$ and $x\vee y >\sqrt t$, $d=2$} 
  
  We continue to suppose  $d=2$ and write $p_t(x)$ for $p_t^{(2)}(x)$ as mentioned previously.
  \begin{prop}\label{thm3.1} \,  Let $d=2$. For each   $M\geq R_A$ and $\e>0$, uniformly for $\x \notin {\rm nbd}_\e(A^r)$ and  $\y $ 
  subject to the constraints $y> \sqrt t$ and $xy <Mt$, as  $y/t \to  0$ and $ t\to\infty$ 
  \beqn\label{16} 
  p^A_t(\x,\y) = \frac{e_A(\x)}{\lg (t/y)}p_t(\y-\x)(1+o(1)).
  \eeqn
  \end{prop}
  \v2\n
 \pf  \,\, The case when  $y/\sqrt t$ is bounded  is covered by Lemmas \ref{lem2.2} and \ref{lem2.3} in conjunction. Hence we suppose that  $y/\sqrt t\to \infty$ as well as $y/t\to 0$ in addition to the constraint  put in the lemma. From $xy<Mt$ 
it then  follows that   $x/\sqrt t \to 0$ and  $x\to \infty$.

Identity (\ref{10}) is valid with $A$ replacing $U(R_A)$, namely
\beqn\label{100} 
p_t(\y-\x) - p_t^A(\x,\y) = \int_0^t ds\int_{\partial A} H_A(\x,s;d\xi)p_{t-s}(\y-\xi).
\eeqn
We are to compute the integral on  the right side  to find the asymptotic form  of $p_t(\y-\x) -p^A_t(\x,\y)$.  Since $P_\x[\sigma_A>t/2] \leq Ce_A(\x)/\lg t$ (see (\ref{17}) below) and
 $p_{t-s}(\y-\xi)$, $\xi\in \partial A$ is $O(p_{t/2}(y))$ for $s>t/2$,
 the integral restricted to   $[t/2,t]$ is  negligible.  Because of the condition   $y/t\to 0$,  $p_{t-s}(\y-\xi)$ may be replaced by $p_{t-s}(y)$.
  It therefore suffices to identify the asymptotic form of 
  \beqn\label{J_012}
  \int_0^{t/2} p_{t-s}(y)P_\x[\sigma_A\in ds].
  \eeqn 
We shall write $J_{[a,b]}$ for this integral restricted to an interval  $[a,b]$, $0\leq a< b\leq t/2$.  Put
 $$\eta=(t/y)^2.$$
  Then $\eta \to\infty$ and $\eta/t\to 0$ in our present setting.  Expanding   $1/(t-s)$ into the Taylor series of  $s/t$ yields
  \beqn\label{taylor}
  p_{t-s}(y) =\frac1{1-\frac{s}t}p_t(y) \exp\bigg\{-\frac{s}{2\eta}\bigg(1 + \frac{s}{t} + \frac{s^2}{t^2}  + \cdots\bigg)\bigg\}.
  \eeqn
Taking  constants $\a>1$  large and $0<\de<1$ small we split the integral   at $\de\eta$ and $\a\eta$ and observe
\[
J_{[0, \de \eta]} =  \int_0^{\de \eta}p_{t-s}(y)P_\x[\sigma_A\in ds] = {p_t(y)} P_\x[ \sigma_A < \de\eta](1+O(\de))
\]
\[
J_{[\a\eta,t/2]} \leq 2p_t(y)\int_{\a\eta}^\infty   e^{-s/2\eta}P_\x[\sigma_A\in ds] 
 \leq  2p_t(y) P_\x[\sigma_A>\a\eta]e^{-\a/2}; \nonumber
 \]
 and
\[\label{J<_>}
J_{[\de\eta,\a\eta]}   \leq 2{p_t(y)} P_\x[\de \eta < \sigma_A < \a\eta]. 
\]
 Let  $\x\notin {\rm nbd}_\e(A^r)$. Then Theorem A.2 entails that  for each  $\a'>0$, uniformly for  $ s > \a' x^2$
 \beqn\label{17}
   P_\x[\sigma_A> s] = \frac{2e_A(\x)}{\lg s}\bigg[ 1 + O\bigg(\frac1{\lg s}\bigg)\bigg].
 \eeqn
Note that  the condition $xy <Mt$ is written as
  $x^2< M^2\eta$ and then  deduce$$\frac{J_{[0, \de \eta]}}{p_t(y)} =  P_\x[ \sigma_A < \de\eta](1+O(\de))  =\bigg(1- \frac{2e_A(\x)}{\lg (\de\eta)}(1+o(1))\bigg)(1+O(\de)),$$
 \[\label{18} 
\frac{ J_{[\a\eta,t/2]} }{p_t(y)}  \leq 2e^{-\a/2}\frac{e_A(\x)}{\lg (\a\eta)}(1+o(1))
 \quad\;\; \mbox{and}\quad\;\;
\frac{ J_{[\de\eta,\a\eta]} }{p_t(y)}  \leq  \frac{4e_A(\x)}{(\lg \eta)^2}[\lg (\a/\de) +O(1)].
\]
 Then, since $1/\de$ and $\a$  may be arbitrarily large, we find that 
     \beqn\label{19}
 J_{[0, t/2]} =\bigg(1- \frac{2e_A(\x)}{\lg \eta}(1+o(1))\bigg)p_t(y).
 \eeqn
Remembering that   the left side equals  $p_t(\y-\x)-p^A_t(\x,\y) - J_{[t/2,t]}$ and  $J_{[t/2,t]}$ is negligible,    we  conclude (\ref{16}), for if  $x^2/\eta \to 0$ (namely, $xy/t \to 0$), then  $p_t(y)\sim p_t(\y-\x)$, whereas
 if  $x^2 > \a' \eta$ for some constant $\a'>0$,  then  $2e_A(\x)/\lg \eta \to 1$, so that $J_{[0,t/2]} = o(p_t(y)) = o(p_t(\y-\x)$). The proof of Proposition \ref{thm3.1} is complete. \qed
 
 \subsection{Case  $xy > t/M$, $d\geq 2$}
 
  \begin{prop}\label{thm2.5} \, Let $d\geq 2$. For each   $M\geq 1$, uniformly for $\x$ and  $\y $ 
  subject to the constraint $xy >t/M$ and $x\vee y <Mt$, as $x\wedge y \to \infty$ 
  $$p^A_t(\x,\y) = p^{(d)}_t(\y-\x) (1+o(1)).$$
  \end{prop}
  
  For the proof we shall use the first half of  the following
  \begin{lem}\label{lem2.1} \, Let $d\geq 1$ and $\nu = \tst12 d-1$. For any $\de>0$ there exists a constant $C_\de$ that depends only on $\de$ and $d$ such that if  $(x+z)z> \de t >0$, then
\v2  
   {\rm (i)} for  $0< \th<1$ 
$$\th^{d/2}\int_{\th t}^t  p_s^{(d)}(x)p^{(d)}_{t-s}(z)ds \leq C_\de \, p_t^{(d)}(x+z)\bigg(\frac{x+z}{tz}\bigg)^\nu\sqrt{\frac{t}{(x+z)z}};
$$

 {\rm (ii)}\, if $z\leq x$ in addition,   $\int_{0}^t  p_s^{(d)}(x)p^{(d)}_{t-s}(z)ds$ is dominated by the right side above.
\end{lem}
\v2\n
\pf\, The computation is based on the identity
\beqn\label{13.0}
p_s^{(d)}(\x)p^{(d)}_{t-s}(\z) =p_t^{(d)}(\z-\x)p^{(d)}_T\Big(\frac{t-s}{t}\x +\frac{s}{t}\z\Big)
\eeqn
where $T= (t-s)s/t$ (recall $p_s^{(d)}(\x)$ is written for  $p_t^{(d)}(x)$).   Our task will be   to  evaluate the integral
$$J:= \int_{\th t}^tp^{(d)}_T\Big(\frac{t-s}{t}\x +\frac{s}{t}\z\Big)ds =\int_{\th t}^tp^{(d)}_T\Big(\frac{t-s}{t}(\x-\z) + \z\Big)ds$$
from above. Choose  the directions of  $\x$ and $\z$ so that $\z$ is a negative multiple of  $\x$; the value is unaltered by this choice, for so is the left side of (\ref{13.0}).  Then   $p_t^{(d)}(\z-\x) = p_t^{(d)}(x+z)$ and  the integrand becomes $ p^{(d)}_T\Big(\frac{t-s}{t}(x+z) -z\Big)$, so that after  substitution of $t-s=u$  
$$J = \int^{(1-\th) t}_0 p^{(d)}_T\Big(\frac{u}{t}(x+z) -z\Big)du, \quad T=\frac{u(t-u)}t.$$
Write the identity  $\int_0^\infty \exp\{-\frac12 \a^2 u - \frac12\b^2 u^{-1}\}u^{-\nu-1}du = 2(\a/\b)^\nu K_\nu(\a\b)$ (cf. \cite{E}, p.146)  as
\beqn\label{eq2.2}
\int_0^\infty p_u^{(d)}(\a u-\b)du =2 (2\pi)^{-d/2}(\a/\b)^\nu K_\nu(\a\b)e^{\a\b},
\eeqn
where $\a\b>0,  \nu = \frac12 d -1$ and $K_\nu$ is the usual modified Bessel function of order $\nu$. Then,  noting  that  $ \th u < T< u$ for $0<u< (1-\th)t$, we have
\begin{eqnarray}\label{00}
 J&\leq& \frac1{\th ^{d/2}}\int_0^\infty p_u^{(d)}\bigg(\frac{x+z}{t}u-z\bigg)du \nonumber \\
  &=& 2\bigg(\frac1{2\pi\th}\bigg)^{d/2}\bigg(\frac{x+z}{tz}\bigg)^\nu K_\nu\bigg(\frac{(x+z)z}{t}\bigg) e^{(x+z)z/t}.
 \end{eqnarray}
Combined  with the following asymptotic formula:
$$K_\nu(\eta)e^\eta \sim \sqrt{\pi/2\eta}\quad \mbox{ as}\quad \eta 
\to \infty$$ 
(for every  $\nu\in \R$)
(cf. \cite{L})  this yields (i). 

Let   $z\leq x$ and   compare the function $ |(t-s)x - sz| $ restricted on $0<s<t/2$ with that on $t/2<s<t$. Then by  symmetry  of $T= s(t-s)/t$ about $s=t/2$   we see that  the integral of $p^{(d)}_T\Big(\frac{t-s}{t}x -\frac{s}{t}z\Big)ds$ on $[t/2,t]$ is not  less than   that on $[0,t/2]$, so  that (ii) follows from (i).
\qed

\v2
{\sc Proof of Proposition \ref{thm2.5}.} We may suppose $A= U(R_A)$. Let  $R_A=1$ and $1< y\leq  x < Mt$. We evaluate the difference $p_t(\y-\x)-p^{U(1)}_t(\x,\y)$ by means of (\ref{10}). As before the integral in (\ref{10}) restricted on $[a,b]\times \partial U(1)$ is denoted by $I_{[a,b]}$.  

 Let $d=2$ and write  $p_t(x)$  for  $p_t^{(2)}(x)$.  Then, using Theorems A.1 and A.4 (of Section 4) we have
\beq
I _{[t/4,t]} &=& \int_{t/4}^t ds\int_{\partial U(1)}  H_{U(1)}(\x,s;d\xi)p_{t-s}(\y-\xi) \\
&\leq& \frac{Ce^{ M}}{1\vee\lg (t/x)} \int_{\partial U(1)} |d\xi|\int_{t/4}^t p_s(\x)p_{t-s}(\y-\xi)ds.
\eeq
 Suppose    $xy>t/M$ and  apply Lemma \ref{lem2.1} with  $\z = \y-\xi$.
Noting  $p_t(\y-\xi-\x) \leq Ce^{2M}p_t(\y-\x)$ we then deduce the bound
\beqn\label{d=2}
\frac{I _{[t/4,t]}}{p_t(\y-\x)} \leq \frac{c_M}{1\vee\lg (t/x)} \sqrt{\frac{t}{xy}}, 
\eeqn
and,  by considering  each case   of   $t/x$ being bounded from below or not, we infer that  the right side above approaches zero as $t\to\infty$.

 For $s\in [0,t/4]$ and $\xi\in \partial U(1)$,  we have $p_{t-s}(\y-\xi) \leq Cp_t(y)$, and  using Theorem A.2 we deduce
 \beq
 I_{[0,t/4]} \leq Cp_t(y)P_\x[\sigma_{U(1)}< t/4]
 \leq  c_Mp_t(y) \frac1{1\vee \lg (t/x)}\cdot \frac{t}{x^2} e^{-2x^2/t}.
 \eeq
 Since $x^2\geq t/M$ and $p_t(y)e^{-2x^2/t} = e^{-y^2/2t}p_t(2x) \leq e^{-x^2/2t}p_t(x+y)$, we  obtain
 \beqn\label{*}
 I_{[0,t/4]} \leq  \frac{Mc_M}{1\vee \lg (t/x)} e^{-x^2/2t}p_t(x+y).
 \eeqn
Thus  $ I_{[0,t/4]}$ is negligible and  for $d=2$, the assertion of the theorem  follows. 

For $d\geq 3$ we apply Theorem A.1 to see that $P_\x[\sigma_{U(1)} <t] \leq C x^{-2}t^{-\nu}p^{(d)}_t(x-1)$ if $x^2\geq t$
and  the same arguments as above yield
$$\frac{I _{[t/4,t]}}{p_t^{(d)}(\y-\x)} \leq \frac{c_M}{y^\nu}\bigg(\frac{x}t\bigg)^\nu\sqrt{\frac{t}{xy}}\quad \mbox{and}\quad  \frac{I_{[0,t/4]}}{p_t^{(d)}(x+y) } \leq  \frac{c_M}{t^\nu} e^{-x^2/2t}$$
in place of (\ref{d=2}) and (\ref{*}), respectively.   The right sides  tending to zero  in  both  inequalities above,      the theorem is thus proved.  
\qed

 \subsection{The higher dimensions  I}
 Let $d\geq 3$ and put
 $$u_A(\x) = P_\x[\sigma_A=\infty].$$
  \begin{thm}\label{thm5.1} \, Let $d\geq 3$. For each $\e>0$, uniformly for $\x,\y \notin {\rm nbd}_\e(A^r)$, as  $t\to\infty$ and $(x\vee y)/t \to 0$
  \beqn\label{thm6_eq}
  p^A_t(\x,\y) = u_A(\x)u_A(\y)p^{(d)}_t(\y-\x) (1+o(1));
  \eeqn
moreover for each   $M\geq 1$,  this equality still holds true  as  $x\wedge y \wedge t \to \infty$ under the constraint  $M^{-1}t< x\vee y< Mt$. 
  \end{thm}

First suppose $x\wedge y\to\infty$. Then,  under the additional restriction  $xy> t$  the assertion of Theorem \ref{thm5.1} follows from   Proposition \ref{thm2.5}. The other case $xy\leq  t$  is easily disposed of by examining the proofs of Lemmas \ref{lem2.2} and \ref{lem2.3}. Indeed in the proof of Lemma \ref{lem2.2} (assuming $x\geq y$ differently from therein)   we obtain  that if $1<y\leq \sqrt t \leq  x\leq t$, then 
$$I_{[0,t/2]} \leq Cx^{-(d-2)} p^{(d)}_{t}(\y-\x)$$
in place of (\ref{101}) since $P_\x[\sigma_{U(1)} <t/2] \leq P_\x[\sigma_{U(1)} <\infty] = x^{-(d-2)}$
and 
$$I_{[t/2,t ]} \leq Cy^{-(d-2)} p^{(d)}_{t}(\y-\x)$$
in place of (\ref{**}) owing  to Theorem A.1  (observe $\int_0^{t/2}p^{(d)}_s(r)ds \leq C r^{-(d-2)}$).
Thus 
$$p^A_t(\x,\y) = p^{(d)}_t(\y-\x) (1+o(1))$$
 if  $x\wedge y \to \infty$ under $xy \leq t$.
The case when $x\wedge y <M $  is now easily  dealt with by looking at (\ref{13}) with  
$r$  suitably chosen  so that $r\to\infty$ subject to $ (x\vee y) r < Mt$,  and arguing as in the proofs
of Lemma \ref{lem2.3} and Theorem \ref{thm2.1}. In view of Proposition \ref{thm2.5} this shows  Theorem \ref{thm5.1}.
 \v2\n
{\sc Remark 2.} In (\ref{thm6_eq})  the error estimate of the order $o(t^{-\eta})$ with some $\eta>0$ in place of $o(1)$ can be derived in a way analogous to the proof of  Theorem \ref{thm2.1}.

 \subsection{Case $y \asymp t$ and $x= O(1)$, $d\geq 2$}

 Here we consider  the case when 
  $x\vee y \asymp t$ and $x\wedge y = O(1)$, which eludes the results stated so far.  
\begin{prop}\label{prop2.1} \, Let $d\geq 2$. For each $\e>0$ and $M>1$, uniformly for    $\x \in \Om_A\setminus {\rm nbd}_\e(A^r)$ and ${\bf v}\in \R^d\setminus\{0\}$ satisfying  $ x<M$ and $\e <|{\bf v}|< M$, as $t\to \infty$ and   $\y/t \to {\bf v}$
$$p_t^A(\x,\y) \sim c^A(\x; {\bf v}) p^{(d)}_t(\y-\x),$$
where $c^A(\x;{\bf v})$ is jointly continuous in  $\x, {\bf v}$, positive and less than unity.
\end{prop}

The next lemma provides  a lower bound, which we need for the proof of Proposition \ref{prop2.1}. It  is 
also used for  the proof of Lemma \ref{lem3.2} of the next section. 
\begin{lem}\label{lem2.9}\, For any $0<\de <1$ and $M>1$, there exists a constant $c_{\de,M}$ such that if  $\de t <R_Ay< M t, x<MR_A$ and $\x \notin {\rm nbd}_{R_A\de}(A^r)$, then
$p^A_t(\x,\y) \geq c_{\de,M} p^{(d)}_t(\y-\x)$ for $t$ large enough.
\end{lem}
\v2\n
\pf\, Suppose $R_A=1$. By the conditions  $\x \notin  {\rm nbd}_{\de}(A^r)$ and  $x<M$ we have  $p^A_s(\xi,\x) >c'=c_{\de, M}'$ for $1/3<s< 2/3$,  $\xi \in U(1+\de)$. Hence, putting  $q_a^{(d)}(y,t) =P_\y[\sigma_{U(a)}\in dt]/dt$, we have for  $t>1$ and $y> 1+\de$
\beq
p^A_t(\x,\y)&\geq& \int_0^t q^{(d)}_{1+\de}(y,t-s)\inf_{\xi\in \partial U(1+\de)}p^A_s(\x-\xi) ds\\
&\geq& c'_{\de, M} \int_{1/3}^{2/3} q^{(d)}_{1+\de}(y,t-s)ds.
\eeq
 But  according to Theorem A.1  we have  $q^{(d)}_{1+\de}(y, t-s) \geq c''_{\de, M}p^{(d)}_t(y)\geq c'''_{\de,M}p^{(d)}_t(\y-\x)$ for $s<2/3$,  $ x<M$ and $\de t < y< M t$  and for $t$ large enough.  Thus the lemma is verified. \qed
 \v2\v2
{\sc Proof of Proposition \ref{prop2.1}.}\, In \cite{Ucalm} it is shown that there exists a measure 
$ \la^A_{{\bf v}}( d\xi)$ on $\partial A$ such that 
$\la^A_{{\bf v}}( d\xi)$ depends on ${\bf v}$ continuously,   the total measure $ \la^A_{{\bf v}}(\partial A)$ is positive, and
\beqn\label{Ucal}
\bigg|\frac{H_A(\y,t;\cdot)}{p_t^{(d)}(y)\La_\nu(y/t)} - \la^A_{{\bf v}}(\cdot)\bigg|_{\rm t.var} \longrightarrow  \, 0
\eeqn
(as $\y/t \to {\bf v}, t\to \infty$), where  $\La_\nu$ is given by (\ref{App}) and   $|\cdots |_{{\rm t.var}}$ designates  the total variation of a signed measure. 
Using this formula we are to compute the integral on the right side of (\ref{100}) with $\x$ and $\y$ interchanged to show  that  for some $c^*= c^*(\x;{\bf v})>0$,
  \beqn\label{cc}
   p_t^{(d)}(\y-\x) -p_t^A(\x,\y) \sim c^*\, p_t^{(d)}(\y-\x). 
   \eeqn
To this end we first see that   this  integral restricted to $(0, t/2]$ is
\[
\int_0^{t/2} ds\int_{\partial A}H_A(\y,s:d\xi)p^{(d)}_{t-s}(\x-\xi)  \leq  \frac{C}{t^{d/2}}P_\y[\sigma_A <t/2], 
\]
hence negligible, for  $P_\y[\sigma_A <t/2] \leq p^{(d)}_{t}(y) \times o(1)$, provided $y\asymp t$.
For the rest of the  integral   we compute $\int_{t/2}^t p_s^{(d)}(\y)p_{t-s}^{(d)}(\x-\xi)ds$ because of (\ref{Ucal}).
Recall  identity (\ref{13.0}) and
 consider the integral  $J$ in the proof of Lemma \ref{lem2.1} with $\th =1/2$ and with $\x$ and $\z$ replaced by $\y$ and  $\x-\xi$ ($\xi\in \partial A$), respectively. On putting  $\tilde \x= \x-\xi$, observe that the integral
 tends to concentrate  on $[t-t^\a,t]$  for any $\a\in (0,1)$ and 
 is asymptotic to $\int_0^\infty p^{(d)}_u\Big( t^{-1}(\y- \tilde \x)u + \tilde \x\Big)du$ (provided that  $y\asymp t$, $x <M$ and $\tilde \x \neq 0$). As in (\ref{eq2.2}) we derive an explicit expression of  this integral, which  immediately leads to
$$ \int_{t/2}^tp^{(d)}_T\Big(\frac{t-s}{t} \y +\frac{s}{t}\tilde \x\Big)ds 
\sim  \frac{2}{(2\pi)^{d/2}}\bigg(\frac{|\y-\tilde \x|}{t| \tilde \x|}\bigg)^\nu K_\nu\bigg(\frac{|\y-\tilde \x||\tilde\x|}{t}\bigg) e^{- (\y-\tilde\x)\cdot (\x-\xi)/t}.$$ 
Now  apply  formula (\ref{Ucal}). Noting $p_t^{(d)}(\tilde \x -\y)= p_t^{(d)}(\y-\x) e^{-\y\cdot\xi/t}(1+O(1/t))$, we then readily  find that   formula (\ref{cc})  holds with 
$$c^*= \frac{e^{-{\bf v}\cdot\x}}{K_\nu(|{\bf v}| )}\int_{\partial A}\frac{K_\nu (|{\bf v}||\x-\xi|)}{|\x-\xi|^\nu} 
 \la^A_{{\bf v}}(d\xi).$$
Thus, in view of Lemma \ref{lem2.9},   the asymptotic formula of the proposition  holds true with  $c^A(\x;{\bf v}) :=1 - c^* \in (0,1)$.   \qed

 \section{Upper and lower bounds in the regime $x\vee y > t$}
\v2\n
Let  $\be =(1,0,\ldots,0) \in \R^d$ and suppose  
$$\y = -y\be  \quad\mbox{and} \quad x\leq y$$ throughout this section. Let  $\x=(x_1, \x')$, where $x_1= \x\cdot \be \in \R$ and $\x' =\pr_{\be}\x =(x_2,\ldots,x_d)\in \R^{d-1}$.
Given $a>0$, put
$$W = W_{a,\be} =\{\z\in \R^d\setminus U(a): \z\cdot \be>0, |\pr_{\be}\z| <a\}$$
and, when  $x_1>0$,  
$$\rho = \frac{x_1}{y+x_1}.$$
Note that $\rho\leq 1/2$.
We shall be  concerned almost exclusively     with the case  $\x \in W$.
We shall decompose $p^A_t(\x,\y)$ by means of the first hitting of  $L_h:= \{(h,\z'): \z'\in \R^{d-1}\}$ ($h\in [0,a]$), the plane perpendicular to the vector  $\be$ and passing through the point $h\be$. 

We are primarily  interested in the case when  $x_1/y <\!\!< 1$ and  $\rho$ can be replaced by $x_1/y$
in   the main results given below. Our choice above  however  is natural  in view of the next lemma;  its proof   may be modified to show that   the first hitting  time of $L_h$ ($0\leq h<x_1$) by the Brownian motion  $(B_t)$ started at $\x$  and conditioned to be in $\y (=-y\be)$ at time $t$  concentrates in  a relatively small interval about $\rho t$ with a  high probability if $\rho$ and $t/(x_1-h)y$ are small enough.   The next  lemma provides   a crude lower bound but covers a wide range of $x_1$ and $y$.
\begin{lem}\label{lem3.1}\, Let $X_t$ be a (standard)  linear Brownian motion and $T_b$ its first passage time of  $b\in \R$.  
  Let  $0<b<\ell $  and  put  $\rho_1 =b/\ell $. Suppose $\rho_1\leq 1/2$.  Then for a universal constant $c>0$,
     \beqn\label{eq9}
P[\, \tst12 \rho_1 t < T_b < \rho_1 t\, |\, X_0=0, X_t=\ell ]  >  c(1\wedge \sqrt{b\ell/t}\,),\quad t>0.
\eeqn
\end{lem}

 \v2\n
\pf \,  Let $k = \ell /t$ and  $T =(t-s)s/t$.  As in the proof of Lemma \ref{lem2.1} the conditional probability  in (\ref{eq9}) is written as  
\beqn\label{eq9a0}
\int^{ \rho_1 t}_{\rho_1 t/2} p^{(1)}_T(k s -b)\frac{b}{s}\,ds.
\eeqn 
Denote  this integral   by  $J$ and observe that 
\beq
J= \int^{\rho_1 t}_{\rho_1 t/2} \exp\bigg\{ -\frac{k b}{2(1-s/t)}\Big(\frac{k}{b} s + \frac{b}{k s}-2\Big)\bigg\}\frac{bs^{-3/2}ds}{\sqrt{2\pi(1- s/t)}}. 
\eeq 
Changing the variable of integration by  $u= (k/b)s = s/\rho_1 t$ we have
\beq
J = \frac{\sqrt{kb}}{\sqrt{2\pi}}\int^{1}_{1/2}
\exp\bigg\{ -\frac{k b}{2(1- \rho_1 u)}\Big(u + \frac{1}{u}-2\Big)\bigg\} \frac{u^{-3/2}du}{\sqrt{1-\rho_1 u}}.
\eeq
Finally, noting $u+u^{-1} -2= (1-u)^2/u \leq 2 (1-u)^2 $ for $u>1/2$,  we infer that  if   $\rho_1<1/2$, then  $J\geq 2 \sqrt{kb} \int_0^{1/2}e^{- kb s^2/(1-\rho_1)}ds\geq 2 \int_0^{\frac12 \sqrt{kb}}e^{-2s^2}ds\geq c(1\wedge \sqrt{kb})$ as desired. \qed

\v2
Denote by $q_0^{(1)}(x,t)$  the density of the distribution $ P[T_x <t\,|\, X_0=0]$, given explicitly by
\beqn\label{FHD}
q_0^{(1)}(x,t) =\frac{x}{t}p_t^{(1)}(x).
\eeqn
 Formula (\ref{eq9}) will be applied in the form
\beqn\label{eq9a}
\int_{\rho_1 t /2}^{\rho_1 t } p_{t-s}^{(1)}(\ell-b) q_0^{(1)}(b,s) ds \geq c\bigg(1\wedge \sqrt{\frac{b\ell}{t}}\,\bigg) p_t^{(1)}(\ell).
\eeqn
Analogously to (\ref{eq9a0}) we have
$$
P[\, \tst12 \rho_1 t < T_b < \rho_1 t\, |\, X_0=0, T_\ell =t ]  = \frac{t}{\ell}\int^{ \rho_1 t}_{\rho_1 t/2} p^{(1)}_T(k s -b)\frac{b(\ell -b)}{s(t-s)}\,ds.$$
Plainly $t(\ell -b)/\ell (t-s) \geq 1-\rho_1$  and  from the  proof  above we obtain 
\begin{lem}\label{lem3.0} \,Under the same setting and assumption as in Lemma \ref{lem3.1}
  \beqn\label{eq9b}
\int_{\rho_1 t/2}^{\rho_1 t }  q_0^{(1)}(\ell-b, t-s) q_0^{(1)}(b, s)ds \geq c\bigg(1\wedge \sqrt{\frac{b\ell}{t}}\,\bigg) \,q_0^{(1)}(\ell, t).
\eeqn
\end{lem}

By the equality $p^{(d)}_{t-s}(\ell-b)p^{(d)}_{s}(b)=p^{(d)}_{t}(\ell)p^{(d)}_{T}(ks-b)$, the  computation  
similar to that carried out for (\ref{eq9a0}) leads to the following 
\begin{lem}\label{lem3.00}  Let  $\a$ be any real constant,  $0<b<\ell $  and suppose  $\rho_1 :=b/\ell \leq 1/2$.   Then for a positive  constant $c>0$ that depends only on $d$ and  $\a$,
$$
\int_{\frac12 \rho_1t}^{\rho_1 t}p^{(d)}_{t-s}(\ell-b)p^{(d)}_{s}(b) s^\a ds\geq c (\rho_1 t)^{\a-\nu}\bigg( \sqrt{\frac{t}{b\ell}} \wedge 1\bigg) \,p^{(d)}_{s}(\ell). 
$$
\end{lem}

These lemmas,  though stated here to explain the role of $\rho$ in the sequel,  are  used  not  for the upper bound given in the next subsection but only for lower bounds.

 Theorem \ref{thm2} follows immediately by combining Propositions \ref{prop3.3} through \ref{prop3.4} given below; $k(\x)$ and $k_*(\x)$ are given by (\ref{450}) and  (\ref{45}), respectively.
\v2
In the sequel  $\k_\de$, $ \k'_\de, \k''_\de$ etc. denote positive constants that depend only on $\de$ and $d$, while  $C, C_1, c, c_1, c'$ etc. continue to denote universal positive constants; they may vary in each  occurrence of them  
\v2\v2
 \subsection{ An upper  bound valid for  $t>0$,   $d=2$}

Put
 \beqn\label{450}
 k(\x) = \frac{(\sqrt{x^2_1 +(a-|\x'|)^2}-x_1)x_1}{(a- |\x'|)^2} \qquad (|\x'| <a, x_1>0).
 \eeqn
It holds that $\sqrt 2-1< k(\x) < 1/2$,  for if $0<\a <1$, then   $\sqrt 2-1< (\sqrt{1+ \a}-1)/\a <1/2$.  
\begin{prop}\label{prop3.3}\, Let $d=2$ and $k(\x)$  be as above. Then there exists a universal constant  $C$ such that   for all  $\x \in W$ and $t>0$, 
$$
\frac{p_t^{U(a)}(\x,\y)}{p^{(2)}_t(\x-\y)} \leq  C\frac{\sqrt {\rho t}}{a-|\x'|}\exp\bigg\{- k(\x)\frac{(a-|\x'|)^2}{\rho t}\bigg\}.
$$
\end{prop}  

\v2\n
\pf\,    We may and do  suppose  $ a-|\x'| \geq  \sqrt{\rho t} $, for if not, the  upper bound to be verified  is trivial in view of  the bound $k(\x) < 1$.

Let $L = \{(0,\xi'): \xi' \in \R\}$, the vertical axis of the plane.  Then, 
\begin{eqnarray}\label{3.11}
p_t^{U(a)}(\x,\y) &\leq& \int_0^t \int_{L \setminus U(a)}p^{(2)}_{t-s}(\y-\xi) P_{\x}[B_{\sigma(L)}\in d\xi, \sigma_L \in ds] \nonumber\\
&\leq& \int_0^t p^{(2)}_{t-s}(\tilde y) P_{\x}[ B_{\sigma(L)}\in L \setminus U(a), \sigma_L \in ds],
 \end{eqnarray}
where $\tilde y = \sqrt{y^2+a^2}$ and we have used $|\y-\xi| \geq \tilde y$ ($\xi\in L\setminus U(a)$) for the second inequality.  Note that
$$
P_{\x}[ B_{\sigma(L)}\in L\setminus U(a),\, \sigma_L \in ds]/ds 
= \frac{x_1}{s}p_s^{(1)}(x_1) P_{\x}[|B'_s|>a], 
$$
where $B_s'$ denotes the second component of $B_s$,
and 
\beqn\label{3.3}
    P_{\x}[\,|B'_s|>a] \leq P_{0} [\, |B'_s|>a- |\x'| \,] \leq \frac{2 s}{a-|\x'|} \cdot \frac{e^{-(a-|\x'|)^2/2s}}{\sqrt{2\pi s}}.
 \eeqn
Then putting  $b=\sqrt{x_1^2+(a-|\x'|)^2}$, we deduce that
\beqn\label{41}
p_t^{U(a)}(\x,\y) \leq  \frac{2x_1}{a-|\x'|} \int_0^t  p^{(2)}_{t-s}(\tilde y)p^{(2)}_s(b)ds.
 \eeqn

Since  $x_1 \geq a- |\x'|$, we have 
$x_1(y+x_1)/t = x_1^2/\rho t  \geq 1$, or what we are interested in, 
$$
 \frac{b(\tilde y+b)}t \geq 1 
$$
by which together with $b\leq \tilde y$ we apply the second half of Lemma \ref{lem2.1} (with $\de=1$) to see that 
\beqn\label{3.3*}
\int_0^{t}  p^{(2)}_{t-s}(\tilde y)p^{(2)}_s(b)ds \leq C_1 p^{(2)}_t(\tilde y+b)  \sqrt{\frac{t}{b(\tilde y+b)}}.
\eeqn
Now, observing   $\sqrt{t/b(\tilde y +b)} \leq C_2\sqrt{\rho t} / x_1$ and
\beqn\label{eqns}
(\tilde y+b)^2 -|\x-\y|^2 \geq 2(b -x_1)y + a^2+ b^2- x^2 \geq 2(b -x_1)(y+x_1) 
\eeqn
  we can conclude that
for some universal constant $C$ 
$$
\frac{p_t^{U(a)}(\x,\y)}{p^{(2)}_t(\x-\y)} \leq  C\frac{\sqrt {\rho t}}{a-|\x'|}e^{-(b-x_1)(y+x_1)/t},
$$
hence the assertion of the proposition since $(b-x_1)(y+x_1)/t =k(\x)(a - |\x'|)^2/\rho t$. \qed


\v2\v2
 \subsection{ Some lemmas in preparation for   lower bounds, $d=2$}
 Let $d=2$.  If $Y_t $ is a  linear Brownian motion,  then for each $ \la\in \R$, we have
\beqn \label{923}
P[Y_{s'}> 0, 0\leq s'\leq  s\,|\,Y_0=\eta_0, Y_s=\eta] = 1-e^{-2\eta_0\eta /s} \quad (\eta_0> 0, \eta>  0, s>0) 
\eeqn
as  is readily derived from the expression of  transition density  for $Y_t$  killed at the origin.
 For $y>0$, $z>0$ and $ t>0$ put
\beqn\label{Q_y}
Q_y(z,t) = q_0^{(1)}(y,t)p_t^{(1)}(z)
\eeqn
($q_0^{(1)}(y,t)$ is given in (\ref{FHD})). From (\ref{923})   it  follows   that for $x=(x_1,x_2)$, $x_1>0$,
\beqn\label{441}
 \frac{P_{\x}[B'_{\sigma(L)} \in dz, \sigma_L\in dt; \forall s\in [0,t], B'_s >0]}{dzdt} =(1-e^{-2x_2z/t} )Q_{x_1}(z-x_2, t).\quad
 \eeqn

 \begin{lem}\label{lem3.2.1}\, For any $0<\de\leq 1$ there exists a constant $\k_\de$ such that for $0\leq \a \leq a$, $y>2a$, $s>\de a^2$ with   $ay/s>1$, and   $z >a+\de\sqrt{as/y}$,
 \beqn\label{eql3.2}
 \frac{P_{-y\be}[B'_{\sigma(L_\a)}\in dz, \sigma_{L_\a}\in ds, \sigma_{U(a)}>s]}{dzds}\geq \k_\de Q_{y+\a}(z,s).\eeqn
 \end{lem}
 \v2\n
 \pf\; The left side of (\ref{eql3.2}) is written as
 $$\int_{-\infty}^\infty  d\eta\int_0^s Q_{y-a}(\eta,s-s')
 \frac{P_{(-a,\eta)}[B'_{\sigma(L_\a)} \in dz, \sigma_{L_\a}\in ds', s'<\sigma_{U(a)}]}{dz}.$$
On restricting the outer integral to the half line $a(1+\de\sqrt{s/ay})\leq \eta < \infty$  and applying  (\ref{441})  this repeated integral is larger  than 
  \beqn\label{lbd_47}\int_{a+\de\sqrt{as/y}}^\infty  d\eta\int_0^s Q_{y-a}(\eta,s-s')(1- e^{-2\de^2as/ys'})Q_{a+\a}(\eta-z, s')ds'
  \eeqn
  if $z >a+\de\sqrt{as/y}$.
Putting $\rho_1= (a+\a)/(y+\a)$  we further restrict   the inner integral to the interval  $0<s' <  \rho_1s $ for which $as/ys' \geq  1/2$; noting  $\rho_1 \leq 2a/(y+a)<2/3$ so that $s'<\frac23 s$ (if $ y>2a$), without difficulty  we also deduce that  for  $z >a+\de\sqrt{as/y}$,
$$\int_{a+\de\sqrt{as/y}}^\infty p^{(1)}_{s-s'}(\eta)p^{(1)}_{s'}(\eta-z)d\eta\geq \k_\de'p^{(1)}_{s}(z),$$
provided $s>\de a^2$ and $\sqrt{as/y} <a$. Hence,  in view of (\ref{Q_y}), the repeated integral in (\ref{lbd_47}) is larger than
$$\k_\de'p_s^{(1)}(z)\int_0^{\rho_1 s}q_0^{(1)}(y-a, s-s')q_0^{(1)}(a+\a, s')ds'.$$
Noting $(y-a)(a+\a)/s  \geq  (y-a)/y > 1/2$ we apply Lemma \ref{lem3.0} to see that the integral above is bounded from below by a positive multiple of $q_0^{(1)}(y+\a, s)$ as desired. \qed
  
  \v2
\begin{lem}\label{lem3.2.2}\, Let  $\tilde B_t=( X_t, Y_t)$ be a standard two-dimensional Brownian motion.   Then,  for all  $k>0$, $\eta >0$, $\eta' > -k\ell$
and $t>0$,
$$P[ Y_s> -k X_s \;\; \mbox{for}\;\; 0<s<t\; |\, \tilde B_0=(0,\eta), \tilde B_t= (\ell, \eta')] = 1- e^{-2\eta (\eta'+k\ell)/(1+k^2)t}.$$
\end{lem}
 \v2\n
 \pf \, From (\ref{923}) we trivially obtain $P[B'_s>0 \; \mbox{for}\; 0<s<t\; |\, B_0= (0, y), B_t = (m, y')] =e^{-2yy'/t}$ for any positive $y, y', m$.  Let $\th\in (0,\pi/2)$ be such that $\tan \th =k$. Then the identity of the lemma is the expression of this one in the coordinate system rotated by $\th$ and shifted by $\eta\sin \th$ along the  horizontal axis, where
  $y= \eta'\cos \th$,
 $y' = (\eta'+k\ell) \cos \th$  and  $m^2 +(y'-y)^2  = l^2+ (\eta'-\eta)^2$.  \qed

\begin{lem}\label{lem3.2.3}\, For  some universal constant $c$ and  for all $\a>0$
and $s>0$,
$$\int_{\a}^{2(\a\vee \sqrt s\,)} p_s^{(1)}(x)dx \geq c\Big(1\wedge \Big[\frac{s}{\a}p_s(\a)\Big]\Big).$$
\end{lem}
\v2\n
\pf\,  Suppose  $\a/\sqrt s \geq  3/2$. Then 
$$\int_{\a}^{2\a} p_s^{(1)}(x)dx \geq \Big[\frac{s}{x}p_s(x)\Big]_{x=\a}^{2\a} - \frac{s^2}{\a^3}p_s(\a)
\geq  \frac59 \frac{s}{\a}p_s(\a) -  \frac{s}{2 \a}p_s(2\a) \ge c\frac{s}{\a}p_s(\a).$$
If  $\a/\sqrt s <  3/2$, then $\int_{\a}^{2 \sqrt s} p_s^{(1)}(x)dx= \int^2_{\a/\sqrt s}p_1^{(1)}(x)dx \geq
 \frac12 p_1^{(1)}(2)$. \qed 
 

\v2\v2
 \subsection{ Lower bounds for the case $\rho t>1$ and proof of Theorem \ref{thm2}}

\v2\v2

\begin{prop}\label{prop3.20}\, Let $d=2$. There exists a universal constant  $c_0>0$ such that  if $\x \in W$, $x_1 >2a$ and  $a^2<t<ay$,  
$$
\frac{p_t^{U(a)}(\x,\y)}{p^{(2)}_t(\x-\y)} \geq  c_0\bigg[\frac{\sqrt {\rho t}}{a-|\x'|} \wedge 1\bigg]\exp\bigg\{- \frac{(a-|\x'|)^2}{2\rho t} [1+ 2a/x_1]\bigg\}.
$$
\end{prop}  

\v2\n
\pf\,    
Let $L_a = \{(a,\xi_2): \xi_2 \in \R\}$.  This time we use the identity  
\begin{equation}\label{3.111}
p_t^{U(a)}(\x,\y) = \int_0^t \int_{-\infty}^\infty p^{U(a)}_{t-s}(\xi, \y) P_{\x}[B'_{\sigma(L_a)}\in d\xi_2, \sigma_{L_a} \in ds] 
 \end{equation}
where   $\xi = ( a,\xi_2)$ ($B'_t$ denotes the second component of $B_t$ as before).  

For the present purpose of obtaining  a lower bound we restrict the outer integral to the interval $\frac12 \rho t \leq s\leq \frac32\rho t$. The strong Markov property applied at the
hitting time   of $L$ yields 
$$
p^{U(a)}_{t-s}(\xi, \y)
 \geq \int_{\frac12 (t-s)a/y}^{(t-s)a/y} d s' \int_{a+ \sqrt{(t-s)a/y}}^\infty Q_{y}(\tilde \xi_2,t-s- s')p^{U(a)}_{s'}(\tilde \xi, \xi)d\tilde \xi_2,
$$
where $\tilde \xi =(0,\tilde \xi_2)$.
 Put $h =\sqrt{(t-s)a/y}$ ($\leq a$).  Since the line passing through the  two points $(0, a+\frac12 h)$ and $(a,a- \frac12 h)$ does not intersect $U(a)$, by Lemma \ref{lem3.2.2} applied with $k= h/a$ ($< 1$) we obtain
$$p^{U(a)}_{s'}(\tilde \xi, \xi) \geq (1-e^{-h^{2}/4s'}) p^{(2)}_{s'}(\xi-\tilde \xi)
\geq c_1p^{(1)}_{ s'}(a)p^{(1)}_{ s'}(\xi_2-\tilde \xi_2)
$$
for $s'<(t-s)a/y$ and, by performing the integration w.r.t.  $\tilde\xi_2$, that if  $\xi_2>a$,
$$p^{U(a)}_{t-s}(\xi, \y) \geq  c_1' p^{(1)}_{t-s}(\xi_2)  \int_{\frac12 (t-s)a/y}^{(t-s)a/y} q_0^{(1)}(y,t-s-s')p^{(1)}_{s'}(a)   d s', $$
provided $s<\frac32 \rho t \, (\leq \frac34t)$. Consequently,  owing to Lemma \ref{lem3.1} (applied with $\de=1/2$)
\beqn\label{eqp_32}
p^{U(a)}_{t-s}(\xi, \y) \geq c_1'' p^{(2)}_{t-s}(\y-\xi)\quad (\xi_2>a). 
\eeqn

We consider separately  the  cases  where $ (a- |\x'|)^2 $ is  less or  not less  than  $\frac32 \rho t $.  First  suppose that  $ a- x_2  \geq  \sqrt{\frac32 \rho t }$ where $x_2=|\x'|$ and restrict  the inner integral  in (\ref{3.111}) to $a\leq \xi_2\leq 2a$. Since  $ |\y-\xi|^2 = (y+a)^2 + \xi_2^2$, this yields for $a< \xi_2<2a$
$$p_t^{U(a)}(\x,\y) \geq c_2\int_{\frac12 \rho t}^{\frac32\rho t}  p^{(2)}_{t-s}(y+a) \frac{x_1-a}{s}p_s^{(1)}(x_1- a) 
P_{\x}[a \leq B'_{s} \leq 2a ] ds,$$
 provided $t>a^2$.
For $ s <\frac32\rho t$ so that $a-x_2> \sqrt{s}$, by Lemma \ref{lem3.2.3}
\beqn\label{intgrl2}
P_{\x}[a \leq B'_{s} \leq 2a ]=  \int_{a-x_2}^{2a-x_2} p^{(1)}_s(\xi_2)d \xi_2
\geq \frac {c'_2 s}{a-x_2}p_s^{(1)}(a-x_2).
\eeqn
 Let $x_1>2a$ and put $\hat y =y+a$, $\hat x_1= x_1-a$ and $b= \sqrt{\hat x_1^2 + (a-x_2)^2}$. Noting $\rho<b/(\hat y+b) \leq \frac32\rho$, we then apply Lemma \ref{lem3.00}  to obtain
 \beqn\label{intgrl}
 p_t^{U(a)}(\x,\y) \geq c_2''\frac{b}{a-x_2} \int_{\frac12 \rho t}^{\frac32\rho t} p_{t-s}^{(2)}(\hat y) p_s^{(2)}(b)ds = c_3\frac{\sqrt{\rho t}}{a-x_2} p_t^{(2)}(\hat y + b).
\eeqn
Now we may proceed  in the same way as in the preceding proof. Indeed, noting
 $|\y-\x|^2 = \hat x_1^2 + \hat y^2 + 2\hat x_1\hat y + x_2^2$ and $b\leq \hat x_1+ (a-x_2)^2/2\hat x_1$,  in place of  (\ref{eqns}) we have
$$(\hat y+b)^2 -|\x-\y|^2 = 2(b -\hat x_1)\hat y + b^2- \hat x_1^2 - x_2^2\leq \frac{(a-x_2)^2}{ x_1 -a}(y+x_1) + O\Big( \frac{1}{x_1^2}\Big),$$
which gives the lower bound of  $e^{-[(a-x_2)^2/2\rho t] x_1/(x_1-a)}$  for $p^{(2)}_t(\hat y+b)/p^{(2)}_t(\x-\y)$.    Finally  using the inequality  $(x_1-a)^{-1} < x_1^{-1}(1 + 2 a x_1^{-1})$  concludes  the required lower bound.

For the case  $0\leq  a- x_2 \leq  \sqrt{\frac32 \rho t} $, the restriction of the inner integral  in (\ref{3.111}) is made to the interval  $a\leq \xi_2 \leq a+4\sqrt {\rho t}$. Applying Lemma \ref{lem3.2.3} again we have
for $\frac12 \rho t\leq s \leq \rho t$
$$P_{\x}[a \leq B'_{s} \leq a +4\sqrt{\rho t} \, ]= \int_{a-x_2}^{a-x_2+ 4\sqrt {\rho t}} p^{(1)}_s(u )d u
\geq c_2$$
 instead of (\ref{intgrl2}), and  a computation using (\ref{eq2.2}) (with $\nu=1/2$)  readily leads to the result. The details are omitted.
\qed

\begin{lem}\label{lem3.2}\, Let $d\geq 2$. For any $\de > 0$, there exists a constant $\k_{\de}$ such that if  $x> (1+\de)a$ and if either $\x \notin W$ or $\x \in W$ with  $\rho t\geq \de a^2 $, then
$p^{U(a)}_t(\x,\y) \geq \k_\de \, p^{(d)}_t(\y-\x)$.
\end{lem}
\v2\n
\pf\, Suppose $\x \in W$ and $x> (1+\de)a$ with  $\rho t> a^2\de$, the case
 $\x \notin W$  being easier (see Proposition \ref{prop3.4} if necessary). Writing  $\x =(x_1, \x')$ (as before)  we can also
 suppose $x_1 >3a$, for otherwise we have $3t\geq \de a y$ so that the result follows from Lemma \ref{lem2.9}.  Let $L_a = \{(a,\z'): \z'\in \R^{d-1}\}$. For simplicity we suppose $\x'=0$, entailing in particular  $p^{(d)}_t(\y-\x)=p^{(d)}_t(y+x_1)$.   Make decomposition
\begin{eqnarray}\label{3.0}
p_t^{U(a)}(\x,\y) = \int_0^t \int_{d\xi\subset L_a}p^{U(a)}_{t-s}(\xi,\y) P_{\x}[B_{\sigma(L_a)}\in d\xi, \sigma_{L_a} \in ds],
 \end{eqnarray}
and note that 
\beqn\label{3.00}
P_{\x}[ B_{\sigma(L_a)}\in d\xi, \sigma_{L_a} \in ds] 
= P[B'_s\in d\xi'\,|\, B_0'= 0]P[T_{x_1-a}\in ds\,|\, X_0=0] 
\eeqn
where  $\xi= (a,\xi')$ and  $T_r$ denotes the first  passage time of  $r$  for a linear Brownian motion $(X_t)$. 
Owing to the assumption  $ \rho t /a^2>\de$, there  exists  a constant $\k'_\de >0$ such that
$$P[2a<| B'_s| < 4 \sqrt t\,|\, B_0'= 0]>\k'_\de \quad\mbox{if}\quad 2^{-1} \rho t < s <2^{-1}3 \rho t. $$
It also holds that if $s< \frac32\rho t \, (< \frac34t)$  and $2a \leq |\xi'|<4\sqrt t$, 
$$p^{U(a)}_{t-s}(\xi,\y) >c\, p^{(d)}_{t-s}(\y-\xi) =c\,  p^{(1)}_{t-s}(y+a)p^{(d-1)}_{t-s}(\xi')\geq c'  p^{(1)}_{t-s}(y+a)p^{(d-1)}_{t}(0)$$ 
for some $c, c'>0$.  Hence  
 \beq
 \frac{p_t^{U(a)}(\x,\y)}{p^{(d)}_t(\y-\x)}
 &\geq& \frac{\k_\de''\, p^{(d-1)}_t(0)} {p^{(d)}_t(y+x_1)}\int_{\rho t/2}^{\frac32\rho t} p^{(1)}_{t-s}(y+a) P[T_{x_1-a}\in ds\,|\, X_0=0] \\
 &=& \k_\de'' P[ \tst12 \rho t < T_{x_1-a} < {\textstyle \frac32}\rho t\,|\, X_0=0, X_t  =  y+x_1] \\
 &>& \k'''_\de,
 \eeq
 where  Lemma \ref{lem3.1} is used  for    the last inequality.
This concludes the proof. \qed

\v2
\begin{prop}\label{prop3.4}\, Let $d=2$. For any $\de>0$  there exists a constant  $\k_{\de}$  such that   if   $\x \in W\setminus U(a+ \de\sqrt{\rho t})$, $t>a^2$  and  $ a-|\x'| \geq \de  \sqrt{\rho t} $, then  
$$
\frac{p^{U(a)}_t(\x,\y)}{p^{(2)}_t(\x-\y)} \geq \frac{\k_\de\,\sqrt {\rho t}}{a-|\x'|}\exp\bigg\{-k_*(\x)\Big(\frac{a-|\x'|}{\sqrt{\rho t}} +\de \Big)^2\bigg\}
$$
with a function  $k_*(\x)$, $\x \in W$ such that 
 $\tst12<k_*(\x)<1$ and $ \lim_{x_1\to\infty}  k_*(\x) = \tst12$. 
\end{prop}

\v2\n
\pf\, The proof is a modification of that of Proposition \ref{prop3.20}. Write $\x=(x_1, x_2)$ as before.  We may suppose $x_2\geq 0$. Let ${\bf u}$ be the  point of intersection where the horizontal  line   $\{(\eta,a): \eta \in \R\}$ meets the half line in the upper half plane that is issuing   from $\x$ and  tangential to  the circle   $\partial U(a)$. Let $u_1$ be the first coordinate 
 of ${\bf u}$ (so that  ${\bf u}=(u_1, a)$) and $L^*$ be the shift of $L:=L_0$ by  $u_1\be$: $L^* =L_{u_1}= \{(u_1,\xi_2): \xi_2\in \R\}$  (see Figure 1).    We decompose  $p^{U(a)}(\x,\y)$ by means of the first hitting of $L^*$: 
  $$
p_t^{U(a)}(\x,\y) = \int_0^t \int_{L^* \setminus U(a)}p^{U(a)}_{s}(\xi,\x) P_{\y}\Big[B_{\sigma(L^*)}\in d\xi, \sigma_{L^*} \in t-ds, \sigma_{U(a)}> t- s\Big].
 $$
(Here the starting site  is not $\x$ but $\y$; this choice is made since our estimation of the hitting distribution of $L^*$ under  $P_\x$ becomes crude  as $\x$ gets close to $a\be$.)
  We restrict the range of  integration  to  the rectangle  $\{(s, \xi_2)\in [0, 4\rho t]\times [a+ \de\sqrt{\rho t}, 2\sqrt t]\}$.  We may and do suppose $\rho \leq1/8$ so that  $4\rho t\leq t/2$, otherwise   $a^2<t< 8a^2/\de^2$ owing to its premise and   the assertion of the lemma being easily shown.

 \begin{figure}[b]
 \begin{center}
\begin{picture}(440,145)(0,40)
\put(236,110){\vector (1,0){154}}
\put(236,130){\line (1,0){154}}
\put(236,110){\line (-1,0){64}}
\put(168,110){\line (-1,0){4}}
\put(160,110){\line (-1,0){4}}
\put(152,110){\line (-1,0){4}}
\put(144,110){\line (-1,0){4}}
\put(136,110){\line (-1,0){4}}
\put(128,110){\line (-1,0){4}}
\put(120,110){\line (-1,0){4}}
\put(112,110){\line (-1,0){4}}
\put(104,110){\line (-1,0){70}}
\put(236,168){\line (0,-1){14}}
\put(236,138){\line (0,-1){58}}
\put(242,168){\line (0,-1){88}}


\put(248,126){\line(3,-2){16}}
\put(248,126){\line(-3,2){20}}
\put(236,110){\line (2,3){11}}

\put(236,110){\circle{40}}
\put(236,110){\circle*{4}}
\put(242,144){\circle*{2}}
\put(242,130){\circle*{3}}
\put(265,115){\circle*{3}}
\put(44,110){\circle*{3}}
\put(44,110){\circle*{3}}

\put(227,113){$0$}
\put(234,143){$\xi$}
\put(393,106){$\be$}
\put(270,114){$\x$}
\put(40,100){$\y$}
\put(244,133){${\bf u}$}

\put(225,160){$L$}
\put(245,160){$L^*$}
\put(15,48){Figure 1: $\bf u$ is  a  point of intersection where the half line  $\{(\eta,a): \eta>0\}$ meets}
\put(64,34){ a  line passing through $\x$ and   tangential to  the circle   $\partial U(a)$.}
\end{picture}
\end{center}
\vspace*{0cm}
\end{figure}


Considering the hitting of an appropriate half line  emanating from the origin  and arguing as in the derivation of  (\ref{eqp_32}) if necessary, we infer that if  $s\leq  4\rho t$ and $\xi\in L^*$ with $\xi_2> a+\de\sqrt {\rho t}$
 $$p^{U(a)}_{s}(\xi,\x) \geq \k_\de' p^{(2)}_{s}(\x -\xi)$$
for some constant $\k_\de'$   .
Hence   
 \beqn\label{3.21}
p_t^{U(a)}(\x,\y) \geq \k_\de' \int_0^{4\rho t}\int_{ a+ \de\sqrt{\rho t}}^{4\sqrt t} p^{(2)}_{s}(\x -\xi) P_{\y}\Big[ B'_{t-s} \in d\xi_2,  \sigma_{L^*} \in t- ds, \sigma_{U(a)}>t-s\Big],
 \eeqn
 where $\xi=(u_1, \xi_2)$.    Owing to  Lemma \ref{lem3.2.1}
  the condition  $\sigma_{U(a)}>t-s$ can be discarded from the  probability in (\ref{3.21})  if  the constant factor is replaced by a smaller one.    Observe that for $s<4\rho t$ and  $0< \xi_2<2\sqrt t$,
 $$P_{\y}[ B'_{t-s} \in d\xi_2,  \sigma_{L^*} \in t- ds] \geq c_1\frac{y}{t} p^{(2)}_{t-s}(y+u_1)d\xi_2ds,$$
  so that  the repeated integral  in (\ref{3.21}) is bounded below by a positive multiple of 
  $$ \frac{y}{t} \int_0^{4\rho t}  p^{(2)}_{t-s}(y+u_1)p^{(1)}_{s}(x_1-u_1) ds \int_{ a+ \de \sqrt{\rho t}}^{4\sqrt t} p^{(1)}_{s}(x_2 -\xi_2)d\xi_2.$$
By  $ a-x_2 \geq \de \sqrt {\rho t}$ and $t>a^2$ it follows that  $a+\de\sqrt {\rho t} <2 \sqrt t$ and  we may apply Lemma \ref{lem3.2.3} to see that for $(a-x_2)^2< s < 4\rho t$,
 \beqn\label{D=2}
 \int_{ a+\de \sqrt{\rho t} }^{4\sqrt t} p^{(1)}_s (x_2 -\xi_2)d\xi_2
\geq \frac{c_2s}{a-x_2+ \de\sqrt{\rho t}} \, p_s^{(1)}(a-x_2+\de \sqrt{\rho t}\,).
\eeqn
Hence, putting 
 $$b^*:= \sqrt{(x_1-u_1)^2 + (a- x_2 + \de\sqrt{\rho t}\,)^2},$$
we have 
  \[
p_t^{U(a)}(\x,\y) \geq  \frac{\k''_\de\,y/t }{a-x_2+ \de\sqrt {\rho t}} \int_0^{4\rho t}  p^{(2)}_{t-s}(y+u_1)p^{(2)}_s(b^*)s ds.
 \]
 
Put $\rho^* =b^* /(y+u_1+ b^*)$.  The condition  $a-x_2> \de\sqrt{\rho t}$ entails $|\x-{\bf u}| < b^*< 2|\x-{\bf u}| $ and  a simple geometric argument leads to
$\frac12 x_1 < |\x-{\bf u}| < 2 x_1$.  Obviously $u_1+b^* >x_1\geq a- x_2$ and these together verify
$$ 2^{-1} x_1 \leq b^* \leq 4 x_1 \quad  \mbox{and} \quad  2^{-1}\rho  < \rho^* <4\rho.
 $$
Hence by Lemma \ref{lem3.00}
$$
\frac{\int_0^{4\rho t}  p^{(2)}_{t-s}(y+u_1)p^{(2)}_s(b^*)s ds}{ p^{(2)}_t(y+u_1+b^*)}\geq c\, \rho^* t\,\sqrt{\frac{t}{(y+u_1+b^*)b^*} \wedge1 }  \geq c' \frac{t}{y}[\sqrt{\rho t} \wedge b^*]
$$
 Since  $b^* \geq \de\sqrt{\rho t}$, we can conclude that
\beqn\label{3.2.5}
p_t^{U(a)}(\x,\y) \geq   \frac{\k_\de \sqrt{\rho t} }{a-x_2+ \de\sqrt {\rho t}} \, p^{(2)}_t(y+u_1+b^*).
\eeqn

Put
\beqn\label{45}
k_*(\x) = \frac{x_1}{x_1-u_1 +\sqrt{(x_1-u_1)^2+(a-x_2)^2}}.
\eeqn
Then   
\beq
(y+u_1+b^*)^2-(y+x_1)^2 = 2y(b^*-x_1+u_1)+ O(1) &=& 2y\frac{(a-x_2+\de\sqrt{\rho t})^2}{x_1-u_1+b^*}+O(1)\\
&\leq&  2k_*(\x)\frac{(a-x_2+ \de\sqrt{\rho t})^2}{\rho},
\eeq
hence 
${p^{(2)}_t(y+u_1+b^*)}/{p_t(y+x_1)}\geq C\exp\{-k_*(\x)(a-x_2 + \de \sqrt{\rho t})^2/\rho t\}$,
which combined with (\ref{3.2.5}) yields the lower bound of $p_t^{U(a)}(\x,\y)$ asserted in the lemma.
Elementary geometry shows that
$(2x_1)^{-1} < u_1 \leq \sqrt{(x_1-u_1)^2+(a-x_2)^2}$ (with the  equality only if $|\x|= a$), hence 
$\frac12<k_*(\x)<1$, finishing  the proof. 
\qed
\v2\v2

 \subsection{ The  higher  dimensions II}
In  higher dimensions  $d\geq 3$ there arises   a  difference in the estimation of the probability
 $$P_\x[ |B_s'|> a]:$$
 for small $s$ this   behaves differently depending on whether $a|\x'|$ is larger than $s$ or not  as is exhibited  in Lemma \ref{lem3.6} below and we accordingly  need to separate  the result into two cases.   Let $k(\x)$  be the function given in  (\ref{450}).

\begin{thm}\label{thm3.5}\, Let  $\x \in W\setminus U(a)$.   There exist  universal constants  $C$ and $c>0$ such that
$$
 \frac{p_t^{U(a)}(\x,\y)}{p^{(d)}_t(\x-\y)} \leq  \left\{ \begin{array} {ll}
{\displaystyle C \Big(\frac{a}{\sqrt {\rho t}\wedge a}\Big)^{d-3}\exp\bigg\{-k(\x)\frac{(a-|\x'|)^2}{\rho t}\bigg\} },
\quad &   a|\x'|< \rho t, \\[5mm]
{\displaystyle  C \frac{a^\nu\sqrt{\rho t}}{|\x'|^\nu(a-|\x'|)}\exp\bigg\{-k(\x)\frac{(a-|\x'|)^2}{\rho t}\bigg\}},
\quad & a|\x'|\geq \rho t;  \end{array} \right.
$$
and that if $x_1 \geq 2a$ and $a^2< t< ay$ in addition,  
$$
 \frac{p_t^{U(a)}(\x,\y)}{p^{(d)}_t(\x-\y)} \geq  \left\{ \begin{array} {ll}
{\displaystyle  c \Big(\frac{a}{\sqrt {\rho t}}\Big)^{d-3}
\exp\bigg\{- \frac{(a-|\x'|)^2}{2\rho t}\Big(1+\frac{2a}{x_1}\Big)\bigg\}},
\quad & a|\x'|< \rho t, \\[5mm]
{\displaystyle c \bigg[\frac{a^\nu\sqrt{\rho t}}{|\x'|^\nu(a-|\x'|)} \wedge 1\bigg]\exp\bigg\{- \frac{(a-|\x'|)^2}{2\rho t}\Big(1+\frac{2a}{x_1}\Big)\bigg\}},
\;& a|\x'|\geq \rho t.  \end{array} \right.
$$
\end{thm}  
\v2
Apparently the upper  and lower bounds given above  correspond to Propositions \ref{prop3.3} and \ref{prop3.20}, respectively. Corresponding result to Proposition \ref{prop3.4} (not stated)  can be derived in the same way as in the proof of Theorem \ref{thm3.5}  given below (see the inequality (\ref{D=2}) for derivation).

\v2
\begin{lem}\label{lem3.6}\, Let $(Y_t)$ be a Bessel process of order $\mu > -1$ and $P^Y_r$  the probability law of $Y$ started at $r\geq 0$. Then for  $0\leq r<a$,
\[
P^Y_r[Y_s > a]  =\left\{\begin{array}{ll} {\displaystyle \frac{a^{2\mu}e^{-(r^2+a^2)/2s}}{2^\mu\Ga(\mu+1) s^\mu}\bigg(1+O\Big(\frac{a^2r^2}{s^2}\vee \frac{s}{a^2} \Big)\bigg)} \quad&\mbox{if} \quad ar< s,\\[5mm]
{\displaystyle\Big(\frac{a}{r}\Big)^{\mu+1/2}\frac{\sqrt s\, e^{-(a-r)^2/2s}}{a-r}
 \bigg(1+O\Big(\frac{s}{(ar)\wedge (a-r)^2}\Big)\bigg) }\quad&\mbox{if} \quad ar\geq  s.
\end{array} \right.
\]
\end{lem}
\v2\n
\pf\, Let $I_\mu$ be the modified Bessel function of the first kind of order $\mu$. Then 
\beqn\label{X}
P^Y_r[Y_s > a]=\int_a^\infty R_s(r,\eta)d\eta,
\eeqn
where for $r>0$,
$$R_s(r,\eta) = \bigg(\frac{\eta}{r}\bigg)^\mu \frac{\eta}{s}I_\mu\bigg(\frac{r\eta}{s}\bigg)e^{-(r^2+\eta^2)/2s}
 $$
and $R_s(0,\eta)= \lim_{r\downarrow 0}R_s(r,\eta)$ (cf. [\cite{RY}). Substitution from the formula
 $$I_\mu(y) =\bigg(\frac{y}{2}\bigg)^\mu\frac1{\Ga(\mu+1)}(1+ O(y^2))\quad (0<y<1); \,\, = \frac{e^y}{\sqrt{2\pi y}}(1+O(1/y))\quad (y>1),$$
yields for $r< a <\eta$, 
\[
R_s(r,\eta) =\left\{\begin{array}{ll} {\displaystyle \frac{\eta^{2\mu +1} e^{-(r^2+\eta^2)/2s}}{2^\mu\Ga(\mu+1)s^{\mu+1}}\bigg(1+O\bigg(\frac{(r\eta)^2}{s^2}\bigg)\bigg)} \quad&\mbox{if} \quad ar< s,\eta <2 a,
\\[4mm]
{\displaystyle\bigg(\frac{\eta}{r}\bigg)^{\mu+1/2}p_s^{(1)}(\eta-r) \bigg(1+O\bigg(\frac{s}{r\eta}\bigg)\bigg) }\quad&\mbox{if} \quad ar\geq  s.
\end{array} \right.
\]
and an  elementary computation of the integral in (\ref{X}) leads to the formulae of the lemma
 (note that if $s>a^2$ (resp. $s>(a-r)^2$)  in the case $ar<s$ (resp. $ar\geq s$), then  the $O$ term becomes dominant so that the  formula  postulated in the lemma is trivial). \qed



 \v2
{\sc Proof of Theorem \ref{thm3.5}.}\,
For both the  upper and lower  bounds  
the proofs of Propositions \ref{prop3.3} and \ref{prop3.4} are readily adapted. Indeed  by examining the proof of Lemma \ref{lem2.1} we can readily obtain
\beqn\label{606}
\int_0^t  p^{(d)}_{t-s}(y)p^{(d)}_s(b)s^\a ds \leq C_{\a,d} p^{(d)}_t( y+b)  \Big(\frac{ y+b}{bt}\Big)^{\nu-\a} \sqrt{\frac{t}{b( y+b)}}.
\eeqn
For the upper bound, in the case $a|\x'|\geq \rho t$
we   apply Lemma \ref{lem3.6} with $\mu= (d-3)/2$ ($=\nu-\frac12$) and the  inequality above with $\a=\nu, b=\sqrt{x_1^2+(a-|\x'|)^2}$ 
(in place of (\ref{3.3}) and (\ref{3.3*}), respectively),  the rest being  virtually unchanged.
As for  the case $a|\x'|<\rho t$ we should take  $\a=0, b=\sqrt{x^2+a^2}$  in (\ref{606})  for which   
$$(y+b)^2- |\y-\x|^2 = a^2 + 2y(b-x_1)\geq \frac{2(a-|\x'|)^2}{\rho} k(\x)- a^2,$$
as is easily checked.
Similar remarks apply  to the lower bound. The details are omitted.

\
 \section{Hitting distributions}
\v2\n

Here we collect several results from \cite{Ubes}, \cite{Ucal_b} and  \cite{Ucalm},  
 the original theorems which they are taken or adapted from  being  indicated in the square brackets.

In the first theorem we put $\nu =\frac12 d-1$ and 
\beqn\label{App}
\La_\nu(y) := \frac{(2\pi)^{\, \nu+1}}{2y^{\nu}K_{\nu}(y)}= \frac{2\pi}{\int_0^\infty  \exp(-\frac1{4\pi u}{y^2 })e^{-\pi u } u^{\nu-1}du}, \quad y>0,
\eeqn
 where $K_\nu$ is the modified Bessel function of the second kind of order $\nu$.
There exists   $\La_\nu(0) := \lim_{y\downarrow 0}\La_\nu(y) \leq \infty$; $\La_0(y) \sim -\pi/\lg y $ \, as $\, y \downarrow 0$, whereas  $0< \La_\nu(0) <\infty$ for $\nu>0$.

\v2\n
{\bf Theorem A.1.} \cite[Theorem 2]{Ubes}.\,  {\it Uniformly for $x > a$, as $t\to\infty$, }
\beq\label{R2}
\frac{P_\x[\sigma_{U(a)} \in dt]}{dt} &\sim &a^{2\nu}\La_\nu\bigg(\frac{a x}{t}\bigg) p^{(d)}_t(x) \bigg[1-\bigg(\frac{a}{x}\bigg)^{2\nu}\,\bigg] \quadd \quad \mbox{if} \quad d \geq 3,\\
&\,&\\
&\sim&p^{(2)}_t(x) 
\times
 \left\{ \begin{array} {ll}  
  {\displaystyle \frac{4\pi\lg(x/a)\,}{(\lg t)^2}   } \quad& (x  \leq \sqrt t\,),\\ [5mm]
 {\displaystyle  \La_0\bigg(\frac{a x}{t}\bigg)   }\quad&( x   > \sqrt t\,),
  \end{array} \right. \qquad \mbox{if} \quad d =2.
\eeq


\v2\n
{\bf Theorem A.2.}
 \cite[Proposition 6.7]{Ucalm}. \, {\it Let $d=2$ and $\e>0$. Then, uniformly for $\x\in \Om_A\setminus {\rm nbd}_\e(A^r)$, as  $t\to \infty$ and $x/t \to 0$}
\beq
 {\rm (i)} \quad  &&P_{\x}[\sigma_A >t] = \frac{2e_A(\x)} {\lg t}\bigg(1+ O\Big(\frac{1}{\lg t}\Big)\bigg)\quad\mbox{ for} \quad x \leq \sqrt t\,;\\
 {\rm (ii)} \quad &&P_{\x}[\sigma_A \leq t] = \frac{1}{2 \lg (t/x)} \int_{x^2/2t}^\infty  e^{-y}y^{-1}dy(1+o(1)) \quad \mbox{for}  \quad  x \geq  \sqrt{t/\lg t}.
 \eeq

  Of the results above Theorem A.1 and  (i) of Theorem A.2 are applied in the proofs of  Propositions \ref{prop2.1} and \ref{thm2.5}, respectively; for all the other applications   we need only the upper  bounds  implied by Theorem A.1 and (ii) of Theorem A.2;  for the latter  only the special case when $A$ is a disc (hence the assertion  is essentially a corollary of Theorem A.1) is relevant.  It is noted that 
  \beqn\label{asymp_e}
  e_A(\x)= \lg (x/R) + E_\x[e_A(B_{\sigma_{U(R)}})]\quad (x\geq R\geq R_A)
  \eeqn
   (cf. \cite[Lemma 6.1]{Ucalm}) and (i) provides a better estimate than (ii) in the range $\sqrt{t/\lg t} \leq x \leq \sqrt t$.

In the following two theorems  $m_a$ denotes the uniform probability measure on $\partial U(a)$.
 \v2\n
{\bf Theorem A.3.} \cite[Proposition 6.1]{Ucalm}. \, {\it  Let $d=2$.  For $r >R_A$ and  $\x\in \Om_A\cap U(r)$,
    \beqn\label{sigma-s}
 P_\x[\sigma_{\partial U(r)} < \sigma_A] =\frac{e_A(\x)}{\lg(r/R_A)}\bigg(1+ \frac{m_{R_A}(e_A)}{\lg(r/R_A)}(1+\de)\bigg)^{-1}
\eeqn   
with $-2(R^{-1}_Ar +1)^{-1}\leq \de \leq 2(R^{-1}_Ar  -1)^{-1}$. Here $m_a(e_A) = \int_{\partial U(a)} e_A(\xi)m_a(d\xi)$.}

 \v2
\v2\n
{\bf Theorem A.4.} \cite[Theorems 2.2 and 2.3]{Ucal_b}. \, {\it   Let   $d\geq 2$.  For  $v\geq 0$,
as  $x/t\to v$ and $t\to\infty$,  
\beqn\label{s-dst}
\frac{P_{\x}[B_{\sigma_{(U(1))}} \in  d\xi\,|\, \sigma_{U(1)}=t]}{m_1(d\xi)}  =
\left\{ \begin{array}  {ll} 1 + O(\ell_d(x,t)x/t)\quad &\mbox{if} \quad v=0,\\[2mm]
 \psi_{v}({\cal R}^\y\xi)(1+o(1))    \quad&\mbox{if}  \quad v>0,
 \end{array}\right.
\eeqn
 uniformly for $(\xi, v) \in  \partial U(1) \times [0, M]$ for each $M>1$.  Here   $\ell_d (x,t) \equiv 1$   for  $d\geq 3$ and}
 $$ \ell_2(x,t) = (\lg t)^2/{\lg(2 \vee x)} \quad (x< \sqrt t); \quad =  \lg(t/x) \quad ( x\geq \sqrt t);$$ 
$\psi_v(\xi)$ is a continuous function of  $(v, \xi)\in [0,\infty)\times \partial U(1)$ that is positive for $v>0$ and $\psi_0\equiv 1$; ${\cal R}^\y$ denotes the rotation that sends $\y/y$ to the unit vector $\be =(1,0,\ldots,0)$ and leaves the plane spanned by $\be$ and $\y$ invariant.
\v2


 \end{document}